%% file: GirstmairKirchner.tex
\documentclass[12pt,a4paper]{article}
\usepackage{amsmath,amssymb}
\def\Bbb{\mathbb}
\def\BQ{``}
\def\EQ{'' }

\title{\bf Towards a Completion of Archimedes' Treatise on Floating Bodies}
\author{Kurt Girstmair and Gerhard Kirchner}
\date{}
\makeatletter
\let\@@maketitle=\maketitle
\def\maketitle{\def\thispagestyle##1{\relax}\@@maketitle}
\makeatother
\newtheorem{theorem}{Theorem}
\newtheorem{proposition}{Proposition}

\def\BE{\begin{equation}}
\def\EDE{\end{equation}}
\def\BD{\begin{displaymath}}
\def\ED{\end{displaymath}}
\def\BA{\begin{array}}
\def\EA{\end{array}}
\def\BEA{\begin{eqnarray}}
\def\EEA{\end{eqnarray}}
\def\BI{\bibitem}

\def\R{\Bbb R}

\def\phi{\varphi}

\def\MB{\mbox}
\def\LD{\ldots}
\def\OV{\overline}

\def\WT{\widetilde}

\def\PP{{\cal P}}
\def\EE{{\cal E}}

\def\BQ{``}
\def\EQ{'' }

\def\MN{\medskip\noindent}

\def\STOP{\hfill$\blacksquare$}

\def\ENS{\enspace}
\def\ET{\WT{E}}
\def\BET{\WT b_1}
\def\BZT{\WT b_2}

\begin{document}
\maketitle

\begin{abstract}
\noindent
\small
In the said treatise Archimedes determines the equilibrium positions
of a floating paraboloid segment, but only in the case when the basis of the segment
is either completely outside of the fluid or completely submerged. Here we give a mathematical
model for the remaining case, i.e., two simple conditions which describe the equilibria in closed
form. We provide tools for finding all equilibria in a reliable way and for the classification
of these equilibria. This paper can be considered as a continuation of Rorres's article \cite{Ro}.
\end{abstract}

%%%%%%%%%%%%%%%%%%%%%%%
\section*{Introduction}
%%%%%%%%%%%%%%%%%%%%%%%

Archimedes' treatise \BQ On Floating Bodies\EQ (its customary Greek title is
$\Pi\epsilon\rho\grave{\iota}$ '$o\chi o\upsilon\mu\acute{\epsilon}\nu o\nu$,
which literally means \BQ about hovering things'', see \cite{Ar}) has been highly esteemed
by mathematicians over
centuries. Book 2 of this treatise can be considered as a sort of crown
of Archimedes' work. In this book he applies a number of his principal results about
volumes and centers of gravity to a problem which is extremely difficult
to handle under Greek premises: namely, the determination of the possible equilibrium positions
of a floating paraboloid segment (for details see Section 1 below).
%throughout the centuries

It seems, however, that even Archimedes was not in a position to treat this
problem in full generality, since he restricts
himself to the cases when the basis circle of the segment either
lies outside of the fluid (i.e., the fluid touches this circle
in at most one point) or is completely submerged.

The case not considered by Archimedes occurs when the basis circle is partially submerged and partially
not. It is, indeed, of a different nature than the \BQ archimedean\EQ case.
Whereas Archimedes gives ruler and compass
constructions for the \BQ tilt angle\EQ of the segment, results of this kind cannot be expected
in the \BQ non-archimedean\EQ case. However, it is possible to establish a {\em mathematical model}
for this case, based on two equations $E=0$ (the equilibrium condition) and $F=0$
(the floating condition). Unlike the corresponding equations in Archimedes'
situation (see Section 1), $E$ and $F$ are no more purely algebraic expressions but also involve the
function $\arctan(x)$. Nevertheless,
these expressions are rather simple if established with care (see Theorem \ref{t1}),
and they do not involve \BQ page-long monstrosities\EQ of which Rorres \cite{Ro} is warning.
% disregarded by Archimedes

The said paper of Rorres contains a graphic completion of the \BQ equilibrium surface\EQ
(based on numerical integration) together with interesting observations of physical phenomena and
inspiring examples (one of which we repeat here in a treatment different of his, see Example 1).
However, \cite{Ro} does not contain any {\em closed formulas} that would describe the non-archimedean case
(which is what we expect from a mathematical model).
From this point of view \cite{Ro} appears only as a {\em first step} towards a completion
of Archimedes' treatise. We hope that the present paper forms a {\em second} step.
Such a step also requires simple tools by which one can decide whether an equilibrium position
is {\em stable} or {\em unstable}. In this connection our above equations
$E=0$ and $F=0$ are, again, quite helpful,
since they imply that the {\em Hesse matrix} of the potential function looks fairly simple for equilibria
(see Section 5).

Maybe the most interesting question in the non-archimedean case concerns
the {\em number of possible equilibria} for a paraboloid segment of a given {\em shape} and
a given (relative) {\em density}.
We have no mathematically rigorous answer
to this question
(which would constitute a {\em third} step).
But we solve a related problem, namely, we determine the number of possible equilibria
if the shape of the segment and the {\em size of the submerged part of the basis} are given,
see Theorem \ref{t2}.
Combined with other devices like Proposition \ref{p2}, this theorem allows
finding \BQ all\EQ possible equilibria for a segment of given shape and density  --- not in
the strict sense of the word but in a convincing manner, as we think.

The standard English translation of Archimedes' treatise seems to be that of Heath \cite{Ar2} from
1897, which is based on Heiberg's first edition of the Greek text.
The German version \cite{Ar3} (with helpful notes) is a translation of
Heiberg's second edition \cite{Ar}, which considers the important Constantinople palimpsest
discovered in 1899.
We highly recommend the monograph \cite{Di} about Archimedes and his works. For a survey
of the treatise on floating bodies the reader may also consult \cite{St}. A number
of problems of floating homogenous bodies are studied in \cite{Cz} and \cite{Gi},
works which also provide a wider theoretical background than we use here.
For additional important references see \cite{Ro}.

%%%%%%%%%%%%%%%%%%%%%%%%%%%%%%%%%
\section*{ 1. The archimedean case}
%%%%%%%%%%%%%%%%%%%%%%%%%%%%%%%%%%%

Throughout this paper we denote subsets of $\R^3$ in an abbreviated way; for instance
\BD
  \{x\le ay+b,\, y\ge c ,\, z\le dx^2\}
\ED
stands for
\BD
  \{(x,y,z)\in \R^3:\,x\le ay+b,\, y\ge c ,\, z\le dx^2\}.
\ED
In this section we give a modern paraphrase of Archimedes' results (see \cite{Ar}, Lib. II) and,
thereby, introduce some basic notations.

It suffices to consider the fixed paraboloid
$
  \{z=x^2+y^2\}
$
that arises from the parabola $\{y=0,\,z=x^2\}$ in the $xz$-plane by rotation around the $z$-axis. Our
{\em paraboloid segment}
$\PP$ is defined by
\BD
  \PP=\{x^2+y^2\le z,\, z\le a\},
\ED
where $a$ is the length of the {\em axis} $\{x=y=0,\,0\le z\le a\}$ of the segment.

Hence the
parameter (in the usual sense) of the rotating parabola equals $1/2$
and the geometric properties of $\PP$ are completely determined
by $a$. The surface of the fluid is a plane
\BD
  \EE=\{z=bx+c\}
\ED
given by the parameters $b,c$. We assume $b\le 0$ throughout this paper.
The segment $\PP$ is said to be in {\em right hand} position if
$\PP\cap \{z> bx+c\}$ lies {\em outside } the fluid (so in Diagram 1 this part of
$\PP$ is on the right hand side of $\EE$). Conversely, for a {\em left hand} position
$\PP\cap \{z< bx+c\}$ must be outside.
In order to exclude uninteresting cases one may also assume that both
$\PP\cap\{z> bx+c\}$ and $\PP\cap\{z< bx+c\}$ are {\em non-empty}.
%contain interior points

%%%%%%%%%%%%%%%%%%
% Diagram 1
\input{Archpict_1}
%%%%%%%%%%%%%%%%%%

In the case considered by Archimedes the intersection of $\EE$ with the {\em basis circle}
$\{x^2+y^2\le a,\, z=a\}$ of the segment $\PP$ consists of {\em at most one point}, which,
in our context,  must be the point $(-\sqrt a,0,a)$.
This situation is henceforth called the {\em archimedean case}.
By our above \BQ non-emptiness\EQ assumption this case excludes a {\em vertical} plane $\EE$, so
our definition of $\EE$ is sufficiently general.

We treat the archimedean case for right hand positions first. To this end we note
the {\em volume} $V$ of $\PP$ and its center of gravity $B$, i.e.,
\BE
\label{1.2}
V=a^2\pi/2 \ENS \MB{ and} \ENS B=(0,0,2a/3).
\EDE
In order to define the  {\em axis} of the paraboloid segment
\BD
  \PP'=\PP\cap\{z\le bx+c\}
\ED
that forms the {\em submerged part}
of $\PP$, we need the midpoint $M=(b/2,0, b^2/2+c)$ of the line segment $\OV{PQ}$,
where $\{P,Q\}$ is the intersection of $\EE$ with the parabola $\{y=0,\,z=x^2\}$.
A vertical line through $M$ intersects this parabola in a point $R$ (the vertex of $\PP'$),
and $R=(b/2,0,b^2/4)$. The axis of $\PP'$ is the line segment
$\OV{MR}$. One immediately finds that its length equals
\BE
\label{1.4}
   a'=b^2/4+c.
\EDE
Archimedes knew that $\PP'$ has the volume $V'=a'^2\pi/2$ (see \cite{Ar}, Lib. II, Sect. iv).
He also knew that its center of gravity
$B'=(x',0,z')$ lies on  $\OV{MR}$ in such a way that $\OV{B'R}$ has length $2a'/3$ (ibid., Sect. ii);
this gives $x'=b/2$, $z'=5b^2/12+2c/3$.

In what follows our fluid has density 1, whereas $\PP$ has density $\sigma$, $0<\sigma<1$
($\sigma=1$ and $\sigma=0$ correspond to the uninteresting cases excluded above).
In this setting Archimedes had to deal with three conditions: First, the {\em floating} condition
$V'=\sigma V$ (this condition is also known as {\em Archimedes' principle}, ibid. Sect. i), which,
because of the above values of $V$ and $V'$,
can be written
\BE
\label{1.6}
 a'=a\sqrt{\sigma}.
\EDE
Second, the {\em equilibrium} condition, which says that $B-B'$ must be perpendicular
to the plane $\EE$. Since we know $(b,0,-1)\perp \EE$, the equilibrium condition
comes down to one of
\BE
\label{1.7}
  b=0 \enspace \MB{ or }\enspace {5b^2}/{12}+{2(c-a)}/{3}+ 1/2 =0.
\EDE
Finally, the condition for the {\em archimedean case}, which reads
\BE
\label{1.8}
 a\ge -b\sqrt a+c.
\EDE
In our formula language these conditions are not difficult to handle.
However, one should be aware of the fact that Archimedes had no formulas at all but could only
work with {\em geometric} propositions, which were enunciated in a rhetorical manner.
In the case $b=0$, (\ref{1.4}) yields $c=a'$, and (\ref{1.6}) gives $c=a\sqrt{\sigma}$.
Since $0<\sigma<1$, we see that $c$ is positive and that (\ref{1.8}) holds automatically.
If $b\ne 0$, (\ref{1.4}) and (\ref{1.6}) yield
\BE
\label{1.5}
  c=a\sqrt{\sigma}-b^2/4,
\EDE
and, thus, (\ref{1.7}) becomes
\BE
\label{1.9}
   b^2=\frac{8a}{3}(1-\sqrt{\sigma})-2.
\EDE
This can only hold if $a >3/4$, so $a\le 3/4$ necessarily requires $b=0$.
Suppose now that $b<0$ is a solution of (\ref{1.9}). As (\ref{1.8}) must be true,
we obtain
$-b\sqrt a\le a(1-\sqrt{\sigma})+b^2/4$
from (\ref{1.5}), and (\ref{1.9}) transforms this inequality into
\BD
    -b\sqrt a\le 5b^2/8+3/4.
\ED
This is true whenever $a\le 15/8$. In the case $a>15/8$, it is the same as saying
that one of
\BD
   0>b\ge \frac{-4\sqrt a}5+\sqrt{\frac{16a}{25}-\frac 65}\enspace \MB{ or } \enspace
   b\le \frac{-4\sqrt a}5-\sqrt{\frac{16a}{25}-\frac 65}
\ED
holds. If we use (\ref{1.9}) in the shape $\sqrt{\sigma}=1-3(b^2+2)/(8a)$, we can read these
inequalities as conditions for $\sqrt{\sigma}$,
namely,
\BE
\label{1.12}
   \sqrt{\sigma}\ge \frac{13}{25}-\frac{3}{10a}+\frac 65\sqrt{\frac{4}{25}-\frac{3}{10a}}
   \enspace \MB{ or }\enspace
   \sqrt{\sigma}\le \frac{13}{25}-\frac{3}{10a}-\frac 65\sqrt{\frac{4}{25}-\frac{3}{10a}}.
\EDE

Summarizing we may say that a {\em right hand equilibrium position} in the archimedean case is possible only
if $a\le 15/8$ or if $\sigma$ satisfies one of the inequalities of (\ref{1.12}). In these cases
either $b=0$ or $b<0$ can be read from (\ref{1.9}), whereas $c$ is given by (\ref{1.5}).
These formulas involve only rational expressions in $a$ and $\sigma$ or square roots of such ones.
Accordingly,  $b$ and $c$ can be constructed by means of ruler and compass if $a$ and $\sigma$ are given.
Archimedes fully described these constructions.
The equilibria defined in this way can be classified as {\em stable} or {\em unstable} by means of the
{\em potential function}, see Section 4.

If we consider $\PP\cap\{z\ge bx+c\}$ as the {\em submerged} part of $\PP$, our right hand position
turns into a left hand one. Since $V-V'$ is the volume of the submerged part,
the floating condition now reads $a'=a\sqrt{\sigma^*}$ for $\sigma^*=1-\sigma$.
In order to obtain an archimedean equilibrium,  $\sigma$ must be replaced by $\sigma^*$ in (\ref{1.9})
and (\ref{1.12}). It is well known that the nature of the equilibrium (stable or unstable)
remains the same, see \cite{Gi}.

%%%%%%%%%%%%%%%%%%%%%%%%%%%%%%%%%%%%%%%%%%%%%%%%%%%%
\section*{2. Equilibria in the non-archimedean case}
%%%%%%%%%%%%%%%%%%%%%%%%%%%%%%%%%%%%%%%%%%%%%%%%%%%%

The case when $\EE$ intersects the basis circle $\{x^2+y^2\le a,\, z=a\}$ of $\PP$ in more than one point
was not considered by Archimedes. Henceforth it will be called the {\em non-archimedean} case.
Note that neither case excludes the other, i.e., for a given paraboloid segment
$\PP$ with given density $\sigma$ archimedean equilibria may occur together with
non-archimedean ones. One can also characterize the non-archimedean case by the condition
\BE
\label{2.1}
   a=bX+c \enspace \MB{ with }\enspace  -\sqrt a < X < \sqrt a.
\EDE
The quantity $X$ is closely connected with the {\em size} of the submerged
part of the basis circle of $\PP$, so it has a natural meaning.
Moreover, for most formulas it is advantageous to use $X$ instead of $c$,
so the reader should get accustomed to the fact that, by (\ref{2.1}),
the plane $\EE$ depends on $b$ and $X$ henceforth, i.e., $\EE=\{z=a+b(x-X)\}$.
Since formulas are slightly simpler for a {\em left hand position}, we assume
that
\BD
 \PP'=\PP\cap \{z\ge bx+c\}
\ED
is the {\em submerged} part of $\PP$. The {\em volume} of $\PP'$ equals $V_1-V_2$, where
$V_1$ and $V_2$ are  the volumes of
\BD
\PP_1=\PP\cap \{x\ge X,\,z\le a\}\ENS \MB{ and }\ENS \PP_2=\PP\cap\{x\ge X,\,z\le bx+c\},
\ED
respectively.
Here $\PP_1$ and $\PP_2$ are no longer paraboloid segments. We call $\PP_1$ a {\em right}
and $\PP_2$ an {\em oblique} paraboloid {\em sector}. It turns out that the volumes of
these sectors are closely related.

In the case of the right sector  $\PP_1$, a cross-section $\PP_1\cap\{x=x_0\}$
is a parabolic segment with area $4(a-x_0^2)^{3/2}/3$. Hence
\BE
\label{2.2}
  V_1=\frac 43\int_X^{\sqrt a}(a-x^2)^{3/2}dx=
  \frac{a^2}2\left(\frac{\pi}2-\arctan \frac{X}{\sqrt A}\right)+\frac{2X^3-5aX}{6}\sqrt{A}
\EDE
with $X$ as in (\ref{2.1}) and
\BE
\label{2.3}
  A=a-X^2>0.
\EDE
As a function of $a$ and $X$, $V_1$ will be denoted by $V(a,X)$.
In order to obtain the
volume $V_2$, we consider another right paraboloid sector, namely,
\BD
   \PP_1'= \{x^2+y^2\le z,\, x\ge X',\,z\le a'\},
\ED
with $a'=b^2/4+c$ defined as in (\ref{1.4}) and
\BE
\label{2.4.1}
X'=X-b/2.
\EDE
By (\ref{2.1}),
$a'$ may be written as
\BE
\label{2.4}
 a'=b^2/4-bX+a.
\EDE
The identity $a'=X'^2+A$ gives $-\sqrt{a'}<X'< \sqrt{a'}$ (so $\PP_1'$ is really
a right paraboloid sector). It is also useful to note
\BE
\label{2.5}
  A=a-X^2=a'-X'^2.
\EDE
Clearly, $\PP_1'$ has the volume $V(a',X')$ in the above sense.
It is not hard to check that the affine mapping
\BE
\label{2.5.1}
  (x,y,z)\mapsto (x,y,z+bx)+(b/2,0,b^2/4)
\EDE
induces a bijection between the right sector $\PP_1'$ and the oblique sector $\PP_2$.
Since the linear part
of this mapping has determinant $1$, $\PP_1'$ and $\PP_2$ have the same volume
\BE
\label{2.6}
V_2=V(a',X')=\frac{a'^2}2\left(\frac{\pi}2-\arctan \frac{X'}{\sqrt A}\right)+\frac{2X'^3-5a'X'}{6}\sqrt A,
\EDE
where we have used (\ref{2.5}). Now the {\em floating condition} (i.e., the
analogue of (\ref{1.6})) reads
\BE
\label{2.7}
 \sigma V=V_1-V_2,
\EDE
with $V=a^2\pi/2$ as in (\ref{1.2}) and $V_1$, $V_2$ as in (\ref{2.2}), (\ref{2.6}), respectively.

Next we consider the center of gravity $B_1=(x_1,0,z_1)$ of $\PP_1$. Its coordinates are defined
by the {\em moments}
\BE
\label{2.8}
   x_1V_1=\frac 43\int_X^{\sqrt a}x(a-x^2)^{3/2}dx, \enspace
   z_1V_1=\frac 4{3}\int_X^{\sqrt a}\frac{3a+2x^2}{5}(a-x^2)^{3/2}dx.
\EDE
This is easy to see, since the above cross-section $\PP_1\cap\{x=x_0\}$ has the center of gravity
$(x_0,0,(3a+2x_0^2)/5)$. Now one can verify (by differentiation of (\ref{2.8}) with respect to
$X$, say) that these moments
take the values
\BE
\label{2.9}
  x_1V_1=\frac{4}{15}A^{5/2}, \ENS z_1V_1=\frac{2a}3V_1+\frac{4X}{45}A^{5/2}.
\EDE
Obviously, these identities hold,
{\em mutatis mutandis}, for the right paraboloid sector $\PP_1'$ as well. In particular,
$\PP_1'$ has the center of gravity $B_1'=(x_1',0,z_1')$, and the analogue of (\ref{2.9})
reads
\BE
\label{2.10}
   x_1'V_2=\frac{4}{15}A^{5/2},\ENS z_1'V_2=\frac{2a'}3V_2+\frac{4X'}{45}A^{5/2}
\EDE
(recall $A=a'-X'^2$, by (\ref{2.5})).

The affine mapping of (\ref{2.5.1}) transforms $B_1'$ into the center of gravity $B_2=(x_2,0,z_2)$
of the oblique paraboloid sector $\PP_2$, so we have
\BD
  x_2=x_1'+b/2,\ENS z_2=z_1'+bx_1'+b^2/4.
\ED
Together with (\ref{2.10}), (\ref{2.4}) and (\ref{2.4.1}) this gives the following formulas for the moments
$x_2V_2$, $z_2V_2$:
\BE
\label{2.12}
  x_2V_2=\frac b2V_2+\frac 4{15}A^{5/2},\ENS
  z_2V_2=\left(\frac{5b^2}{12}-\frac{2bX}3+\frac{2a}3\right)V_2+\left(\frac{4X}{45}+\frac{2b}9\right)A^{5/2}.
\EDE
Let $(x',0,z')$ be the center of gravity of $\PP'$. Obviously, our moments satisfy
\BE
\label{2.13}
 x'(V_1-V_2)=x_1V_1-x_2V_2, \ENS z'(V_1-V_2)=z_1V_1-z_2V_2.
\EDE
We are now in a position to enunciate the analogue of (\ref{1.7})
in the non-archimedean case.
As $B-B'$ must be perpendicular to $\EE$, we obtain
\BD
   2a/3-z'=x'/b.
\ED
Since the volume $V_1-V_2$ of $\PP'$ is positive, this is equivalent to
\BE
\label{2.14}
  (2a/3-z')(V_1-V_2)-(x'/b)(V_1-V_2)=0.
\EDE
Here we insert the right hand sides of (\ref{2.13}) for $x'(V_1-V_2)$ and $z'(V_1-V_2)$.
Moreover, we use (\ref{2.9}) and (\ref{2.12}). Then a short calculation shows that
the left side of (\ref{2.14})
comes down to
\BE
\label{2.17}
  E=\left(\frac{5b^2}{12}-\frac{2bX}{3}+\frac 12\right)V_2+\frac{2b}9 A^{5/2}
\EDE
(here $E$ stands for {\em equilibrium} ).
We recall (\ref{2.7}) and summarize our results in

\begin{theorem} % Theorem 1 %%%%%%%%%%%%%%%%%%%%%%%%%%%%%%%%%%%%%%%%%%%%%%%%
\label{t1}
In the non-archimedean case a left hand equilibrium position is characterized by
the conditions $F=0$ and $E=0$,
where
\BD
F=V_1-V_2 -\sigma V
\ED
and $E$ is given by {\rm(\ref{2.17})}.
The quantities $V$, $V_1$, $V_2$ and $A$ are as in {\rm
(\ref{1.2}), (\ref{2.2}), (\ref{2.6})} and {\rm(\ref{2.3})}.

\end{theorem} %%%%%%%%%%%%%%%%%%%%%%%%%%%%%%%%%%%%%%%%%%%%%%%%%%%%%%%%%%%%%%

\MN
Remarkably, the equilibrium condition $E=0$ involves only the volume $V_2$,
whereas the floating condition $F=0$ also involves  $V_1$. Further, the
expression
\BE
\label{2.18}
f={5b^2}/{12}-{2bX}/{3}+1/2
\EDE
occurring in $E$ is, by virtue of (\ref{2.1}), identical with the expression
$5b^2/12+2(c-a)/3+1/2$ from the equilibrium condition
(\ref{1.7}) in the archimedean case.
It seems, however, that this formal analogy has no influence on the rather
different properties of both cases.

For a {\em right hand} position of $\PP$
and density $\sigma$, the floating condition reads $V_1-V_2-\sigma^*V=0$
with $\sigma^*=1-\sigma$ (see the Archimedean case).
The equilibrium condition remains unchanged.

%%%%%%%%%%%%%%%%%%%%%%%%%%%%%%%%%%%%%%%%%%%%%%%%%%%%%%%%%%%%%%%%%%%%%%%%%%%
\section*{3. Finding \BQ all\EQ solutions in the non-archimedean case}
%%%%%%%%%%%%%%%%%%%%%%%%%%%%%%%%%%%%%%%%%%%%%%%%%%%%%%%%%%%%%%%%%%%%%%%%%%%

As above, suppose $a$ and $\sigma$ are given. Then we know that $\PP$ can take at most two
archimedean right hand equilibrium positions, corresponding to $b=0$
or $b=-(8a(1-\sqrt{\sigma})/3-2)^{1/2}$, see (\ref{1.7}), (\ref{1.9}).
In the non-archimedean situation the determination of the number of possible equilibria
in a mathematically rigorous way seems to be much more difficult, and even a
general upper bound for this number is out of reach for us. In fact, we do not know
how to bound the number of zeros of the (non-algebraic) map
\BE
\label{3.1}
 ]-\infty,0[\:\times\: ]-\sqrt a,\sqrt a[\to \R^2:(b,X)\mapsto (F,E).
\EDE
%which is not an algebraic one.

If, however, we assume that $X$ (instead of $\sigma$) is given together with $a$,
we can determine the exact number of values
$b$ satisfying $E=0$ quite well (see Theorem \ref{t2}).
Each of these values $b$ gives, because of $F=0$,
exactly one density $\sigma$.
In this way we obtain reliable diagrams connecting
$X$ and $\sigma$. For this reason a {\em reliable} (though not rigorous) answer to the question
about the number of zeros of (\ref{3.1}) seems to be possible.
As we remarked already,
$X$ can stand for the size of the submerged
part of the basis circle of $\PP$;
hence replacing $\sigma$ by $X$ is not quite unnatural but has a certain value of its own.

Recall that $E=f\,V_2+2bA^{5/2}/9$ with $f$ as in (\ref{2.18}). Put
\BD
   \ET=\frac{E}{f\,a'^2}=
   \frac{1}2\left(\frac{\pi}2-\arctan \frac{X'}{\sqrt A}\right)+\frac{2X'^3-5a'X'}{6a'^2}\sqrt A+
   \frac{2bA^{5/2}}{9f\,a'^ 2},
\ED
with $a'$, $X'$, $A$ as in (\ref{2.4}), (\ref{2.4.1}), (\ref{2.3}), respectively.
The advantage of $\ET$ lies in the fact that its derivative with respect to $b$ is a {\em rational}
function of $b$. Indeed,
\BD
   \frac{\partial\ET}{\partial\, b\:}=P\,\frac{A^{5/2}}{108\,a'^3f^2},
\ED
where $P$ is the polynomial
\BE
\label{3.3}
 P=6X\,b^3+(-10a+21)\,b^2-36X\,b +12a+18,
\EDE
which is cubic in $b$.
We further note that $f$ has zeros $b_1\le b_2<0$ only if $X\le -\sqrt{15/8}$.
In this case these zeros read
\BE
\label{3.4}
  b_1=\frac{4X}5-\frac 15\sqrt{16X^2-30}, \ENS b_2=\frac{4X}5+\frac 15 \sqrt{16X^2-30}.
\EDE
With these tools at hand, we are able to enunciate the main result of this section,
which describes the (negative) solutions $b$ of $E=0$.

\begin{theorem} % Theorem 2 %%%%%%%%%%%%%%%%%%%%%%%%%%%%%%%%%%%%%%%%%%%%%%%%
\label{t2}
Let $a>0$ and $X$ be given, $-\sqrt a<X<\sqrt a$.

 {\rm (a)} If $X\le -\sqrt{15/8}$, $P$ has exactly two negative zeros $\BET<\BZT$.
Moreover, $\BET<b_1$, and $E=0$ has exactly two solutions,
one in $]\BET,b_1[$, the other in $]b_2,0[$, where $b_1$,
$b_2$ are as in {\rm(\ref{3.4})}.

 {\rm (b)} In the case $-\sqrt{15/8}<X<0$ the equation $E=0$ has solutions only if $P$ has two negative zeros
$\BET<\BZT$. If $E(\,\BZT)=0$, then $\BZT$ is the only solution of $E=0$.
If $E(\,\BZT)<0$, then $E=0$ has exactly
two solutions, which lie in the intervals $]\BET,\BZT[$ and $]\BZT,0[$. If $E(\,\BZT)>0$,
$E=0$ has no solution.

 {\rm (c)} In the case $X=0$ the equation $E=0$ has no solution for $a\le 21/10$.
If $a>21/10$, $E=0$ has exactly one solution,
which lies in $]\BET,0[$, where $\BET$ is the negative zero of $P$.

{\rm  (d)} In the case $X>0$, $P$ has exactly one negative zero $\BET$ and $E=0$ has exactly one solution,
which lies in $]\BET,0[$.

\end{theorem} %%%%%%%%%%%%%%%%%%%%%%%%%%%%%%%%%%%%%%%%%%%%%%%%%%%%%%%%%%%%%%

\MN
Our {\em Proof} of Theorem \ref{t2} requires the knowledge of the limits of $E$
for $b$ tending to $0$ or to $-\infty$. We note

\begin{proposition} % Proposition 1 %%%%%%%%%%%%%%%%%%%%%%%%%%%%%%%%%%%%%%%%%%%%%%%%%%%
\label{p1}

In the above setting,
\BD
\lim_{b\to 0}E= \frac{V_1}2>0 \ENS\MB{ and }\ENS \lim_{b\to-\infty}E=\frac{-4X}{45}A^{5/2}.
\ED
\end{proposition} %%%%%%%%%%%%%%%%%%%%%%%%%%%%%%%%%%%%%%%%%%%%%%%%%%%%%%%%%%%%%%%%%%%%%

\MN
{\em Proof.} The value of the first limit is clear from (\ref{2.17}) and the fact that $V_2$
tends to $V_1$ for $b\to 0$. As to the second one, recall that $\PP_2$ has the center of gravity
$(x_2,0,z_2)$, where $x_2$, $z_2$ satisfy (\ref{2.12}).
Now the second identity of (\ref{2.12}) can be written
\BD
  z_2V_2=E+\left(\frac{2a}3-\frac 12\right)V_2+\frac{4X}{45}A^{5/2}.
\ED
Since $V_2\to 0$ for $b\to-\infty$ and $0\le z_2\le a$, we conclude that $z_2V_2$
tends to zero and $E$ to the value in question.
\STOP

\MN
{\em Proof of Theorem \ref{t2}}. Note $a'>0$ for all $X$ in question, so only the zeros $b_1\le b_2$
of $f$ (as given in (\ref{3.4})) can be poles of $\ET$.

This situation occurs in case (a), which we discuss first. The definition of $\ET$ shows
\BD
%\label{3.10}
   \lim_{b\to-\infty}\ET=0, \ENS \lim_{b\to b_1,\, b<b_1}\ET=-\infty.
\ED
By Proposition \ref{p1}, $\ET$ must be positive for all $b\ll 0$; hence it takes a positive
maximum in $]-\infty,b_1[$. This requires that ${\partial\ET}/{\partial\, b}$ and, thus, $P$,
has a zero $\BET<b_1$, for which the said maximum is taken. By Descartes' rule (see \cite{Ko},
p. 310), $P$ has exactly one positive zero, so it must have another zero $\BZT<0$. We may assume
$\BET\le \BZT$. In the case $X=-\sqrt{15/8}$ one verifies $\BZT=b_1=b_2$. If $X<-\sqrt{15/8}$,
we have $b_1<b_2$ and
\BD
  \lim_{b\to b_1,\, b>b_1}\ET\:=\,\infty\,= \lim_{b\to b_2,\, b<b_2}\ET;
\ED
from this we conclude that $\ET$ takes a minimum in $]b_1,b_2[$, more precisely,
for $b=\BZT\in]b_1,b_2[$.
In this way we know the intervals where $\ET$ is strictly monotonous.
We see, first,
that $\ET$ has exactly one zero $<b_1$, which lies in $]\BET,b_1[$. Since
\BD
  \lim_{b\to b_2,\, b>b_2}\ET=-\infty, \ENS\lim_{b\to 0}\ET>0
\ED
(recall Proposition \ref{p1}),
we see, second, that $\ET$ has exactly one zero in $]b_2,0[$.
As $\ET$ is positive in $]b_1,b_2[$, our list of (negative) zeros of $\ET$ is complete.

In the remaining cases poles of $\ET$ occur no longer. In case (b) the polynomial $P$
also has exactly one positive zero. If $P$ has no zeros $<0$, ${\partial\ET}/{\partial\, b}$ is positive
in $]-\infty,0[$, so $\ET$ is strictly increasing. Since both limits of Proposition \ref{p1} are
positive, $\ET$ does not vanish in $]-\infty,0[$. This assertion remains true if $P$ has a double zero $<0$.
Suppose, therefore, $P(\BET)=0=P(\BZT)$ for $\BET<\BZT<0$. It is not hard to see that $\ET$ must have a
positive maximum at $\BET$ and a minimum at $\BZT$.
Monotonicity arguments and consideration of the sign of $E(\BZT)$ yield the
number and the location of the solutions of  $E=0$ just as indicated in the theorem.

In case (c), $P$ becomes a quadratic polynomial, which is positive for all $b<0$ as long as $a>21/10$.
In this case $\ET$ is strictly increasing in $]-\infty,0[$, and since it tends to zero for
$b\to -\infty$, it must be positive throughout. If $a<21/10$, $P$ has exactly one zero
$\BET$ in $]-\infty,0[$. Now $\ET$ is decreasing for $b<\BET$ but
increasing for $b>\BET$, so it takes a minimum for $b=\BET$. The limit of $\ET$ for $b\to-\infty$
being $0$, we conclude $\ET(\BET)<0$. These arguments imply that
$\ET$ has exactly one zero $<0$, which lies in $]\BET,0[$.

In case (d), Descartes' rule, when applied to the polynomial $P(-b)$,
shows that $P$ has exactly one zero
$\BET<0$. Since $E<0$ for all $b\ll 0$ and $\lim_{b\to-\infty}\ET=0$, the function $\ET$
has a negative minimum at $\BET$. Therefore, $\ET$ vanishes for
some $b\in]\BET,0[$, and this $b$ is the only zero of $\ET$
for reasons of monotonicity.
\STOP

\MN{\em Example 1.} We apply Theorem \ref{t2} in the case considered by Rorres \cite{Ro}: He took
$\sigma=0.51$ and a \BQ base angle\EQ $\phi=74.33$. This value of $\phi$ corresponds to
$a=(\tan\phi)^2/4\approx 3.17690918$ in our setting. From Theorem \ref{t2} we know that, for each
$X\in]-\sqrt a,-\sqrt{15/8}\,]$, the equation $E=0$ has exactly two solutions $b\in]-\infty,0[$.
Hence we expect that $E=0$ describes a curve $(X,b)$ with two branches for all $X$ in this range
(note $-\sqrt a\approx -1.782$, $-\sqrt{15/8}\approx-1.369$). Next our theorem
suggests testing whether $E(\BZT)<0$ for
a sufficiently large number of (equidistant) values of $X$ in $]-\sqrt{15/8},0[$.
Since this is true in all cases, the curve $(X,b)$ will also have two branches for $X$ in
this interval. Finally, the theorem implies that $(X,b)$ consists of one branch for $X\in[0,\sqrt a[$.
Using the floating condition $F=0$ of Theorem \ref{t1} in the shape $\sigma=(V_1-V_2)/V$,
each pair $(X,b)$ produces a pair $(X,\sigma)$. One expects that each branch of $(X,b)$ produces
a branch of the curve $(X,\sigma)$ in this way.

In order to obtain Diagram 2, one has to investigate all
values $X$ between $-1.78$ and $1.78$ in steps of $1/100$.
Numerical values for the corresponding solutions $b$ of $E=0$ can be found by standard methods
as Newton's algorithm.
The diagram displays the resulting points $(X,\sigma)$.
Our steps are small enough to produce the picture of a nearly continuous (and smooth) curve --- in contrast
to the situation of Example 2 below.
The horizontal line $\sigma=0.51$ intersects
the curve in {\em one} point of the upper branch
and in {\em at most two} points of the lower one.
Hence we expect at most three left hand equilibrium positions to be detected in this way.

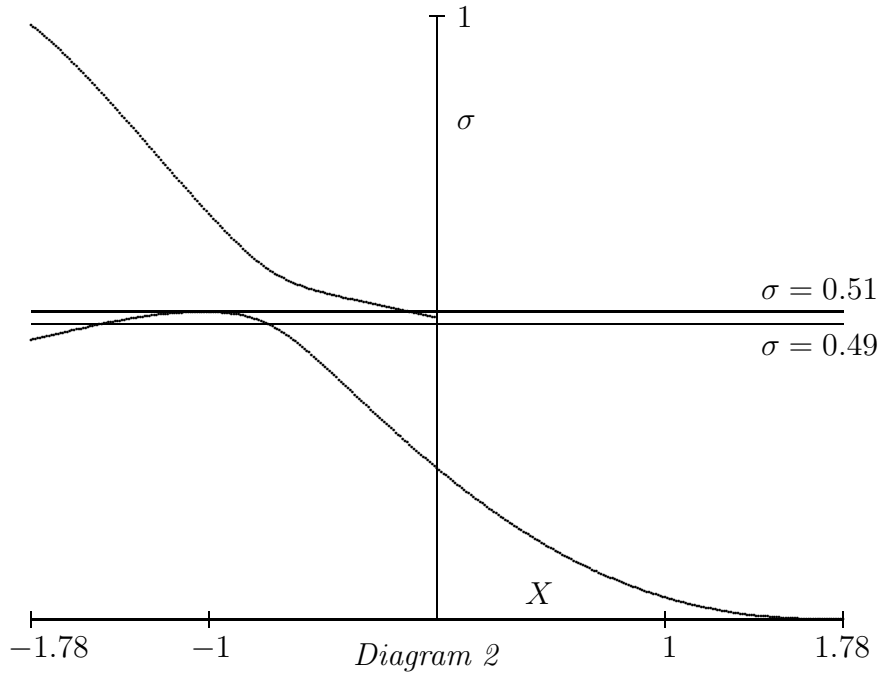
\input{Archpict_2}

%\unitlength1cm
%\begin{picture}(12,10) %Breite und Höhe

% Koordinaten wie üblich: (11,1) heißt im Bild 11(cm) nach rechts und 1(cm)
% nach oben, bei Geraden ist (a,b) der Anstieg b/a und {c} die Länge,

%\put(1,1){\line(1,0){10}}
%\put(1,1){\line(0,1){8.5}}
%\put(11,1){\line(0,1){8.5}}
%\put(1,9.5){\line(1,0){10}}
%\put(4,0.5){\hspace{1cm} {\em Diagram 3}}

%\end{picture}

Inspecting our computation more closely we find the pairs
\BD
(X,\sigma)\approx (-1.04, 0.50997999),\:(-1.03, 0.51000418),\:(-1.02, 0.50999785)
\ED
on the lower branch. Accordingly, this branch should contain
{\em two} equilibria with $\sigma=0.51$. With these starting values and two more
on the upper branch it is no more difficult to
find
\begin{eqnarray*}
  (X,b)&\approx&(-1.03304236, -1.12424322), (-1.02105684, -1.13986072)\\
  \MB{ and }&& (-0.12106085,-12.68795681)
\end{eqnarray*}
(these numerical values of $a$, $X$ and $b$ satisfy $E=0$ and $\sigma=0.51$ up to
an error $<10^{-8}$). The
corresponding approximate \BQ tilt angles\EQ (in Rorres's sense) are
$131.653^{\circ}$, $131.260^{\circ}$ and $94.506^{\circ}$, respectively.
For right hand equilibria
we have to replace $\sigma$ by $\sigma^*=1-\sigma$, which means that we intersect our curve with the
horizontal line $\sigma=0.49$. In this way we detect two more solutions
\BD
(X,b)\approx(-1.46372405,-0.69920557),\ENS (-0.74316119, -1.52773443),
\ED
which
correspond to approximate tilt angles of $34.961^{\circ}$ and $56.793^{\circ}$.
Altogether, we have found five
non-archimedean equilibria in this case; Rorres has only four, since he considers
those corresponding to angles of $131.653^{\circ}$ and $131.260^{\circ}$ as only one solution
with an angle of $131.5^{\circ}$.
Rorres' situation can be established by a slight modification of $\sigma$.
Indeed, $\sigma\approx 0.51000554$ melts those two solutions into
$(X,b)\approx (-1.02702703,-1.13205421)$, the corresponding tilt angle being $\approx 131.456^{\circ}$.
However, Rorres writes $\sigma=0.510$ and $\phi=74.330$, which suggests interpreting
these figures as {\em exact} values.

\MN
{\em Remarks.} 1. The reader need not worry about the fact that a change of the last digit of our values
for $X$ or $b$ may give slightly better results in combination with the above numerical approximation
$3.17690918$ of $a$.
Our results have been obtained by rounding off the higher digits of substantially better values
--- not only of $X$ and $b$ but also of $a$; accordingly,
$X$ and $b$ may appear not sufficiently precise when they are
combined with {\em this} approximation of $a$.

2. Computations suggest that the number of five non-archimedean equilibria
(as in the example) could represent the maximum.

\MN
In this example {\em each} value $X\in]-\sqrt a,\sqrt a\,[$ defines at least {\em one} point
$(X,b)$ that satisfies the equilibrium condition $E=0$. This is no longer true if $a\le 3$.
In this case an investigation of the {\em discriminant} of the polynomial $P$ of (\ref{3.3}) exhibits
values $a>0$ and $X<0$ for which $\ET$ is monotonically increasing; by Proposition \ref{p1}, this means
that $E=0$ has no solution $b\in]-\infty,0[$. We summarize these observations in
%We do not discuss the details of this investigation but only note the result.
%which may be helpful when studying concrete examples:

\begin{proposition} % Proposition 2 %%%%%%%%%%%%%%%%%%%%%%%%%%%%%%%%%%%%%%%%%%%%%%%%%%%
\label{p2}

Let $a>0$ and put $a_1=(-213+198\sqrt 11)/250$ ($\approx 1.7748$). Then the equilibrium
condition $E=0$
%(see (\ref{2.17}))
has no solution $b\in]-\infty,0[$ if

   {\rm (a)} $a\le a_1$ and $X\in]-\sqrt a, 0[$,

   {\rm (b)} $a\in ]a_1, 21/10]$ and $X\in[{X_1},0[$ or

   {\rm (c)} $a\in ]21/10,3]$ and $X\in[{X_1},{X_2}\,]$,
\\ where $X_1=-(\gamma+\sqrt{\delta}/27)^{1/2}$, $X_2=-(\gamma-\sqrt{\delta}/27)^{1/2}$ with
$\gamma=-11a^2/54+5a/9+13/24$ and $\delta=(3-a)(a+6)^3$.

\end{proposition} %%%%%%%%%%%%%%%%%%%%%%%%%%%%%%%%%%%%%%%%%%%%%%%%%%%%%%%%%%%%%%%%%%%%%

\MN
{\em Proof.} We briefly sketch the main arguments. Suppose $X<0$.
Then $P$ has a zero $>0$. Further,
$\ET$ is monotonically increasing (on $]-\infty, 0[$) if, and only if,
$P$ has either no zero $<0$ or a double zero $<0$. This is
the same as saying that the discriminant $D$ of $P$ is $\le 0$.
But $D$ can be written as a quadratic polynomial $D_1$
in $Y=X^2$. The complex zeros of $D_1$ are $Y_1=\gamma+\sqrt{\delta}/27$ and $Y_2=\gamma-\sqrt{\delta}/27$.
These zeros are real only if $a\le 3$, which means that only in this case $D\le 0$ is possible.
Assuming $a\le 3$ henceforth, one sees that $Y_1$ is always positive,
whereas $Y_2$ is positive just if $a>21/10$.
This gives the intervals for $X$ in the cases (b) and (c). However, one has also to observe
that $Y_1< a$ if, and only if, $a> a_1$ (whereas $Y_1\le 15/8$ is always true).
\STOP

\MN
{\em Example 2.} Consider $a=5/2$, which falls under case (c) of the proposition.
Accordingly, the curve $(X,b)$ is empty for $X\in[X_1,X_2]$,
where $X_1\approx -1.143$ and $X_2\approx-0.0917$.
Since $\sqrt{5/2}\approx1.581$ and $\sqrt{15/8}\approx 1.369$, we know from
Theorem \ref{t2} that it should have two branches for $-1.58\le X \le-1.37$.
By means of this theorem one verifies
that these two branches should extend as far as $X=-1.29$.
Furthermore, two branches are to be expected for $-0.08<X<0$, but only one
branch for $X\ge 0$.

The diagram displays points of the curve $(X,\sigma)$ with $-1.58\le X\le 1.58$,
again in steps of $1/100$.
A comparison with Example 1 reveals some marked differences.
First, the two branches for $-1.58\le X\le -1.29$ are {\em connected}, and the same holds
for $-0.08\le X<0$. Second, the diagram
suggests that there is exactly {\em one}  non-archimedean solution for $0.415<\sigma<0.585$
(observe that possible right hand equilibria with $-1/2<\sigma< 0.585$ correspond to
left hand ones with $0.415<\sigma<1/2$; for $\sigma=1/2$ see Section 6).
The density of points is rather low when the tangent of the curve is nearly
vertical, so a \BQ continuous\EQ picture requires much smaller steps for values of $X$
in this region
(say steps of $1/10000$ instead of $1/100$).

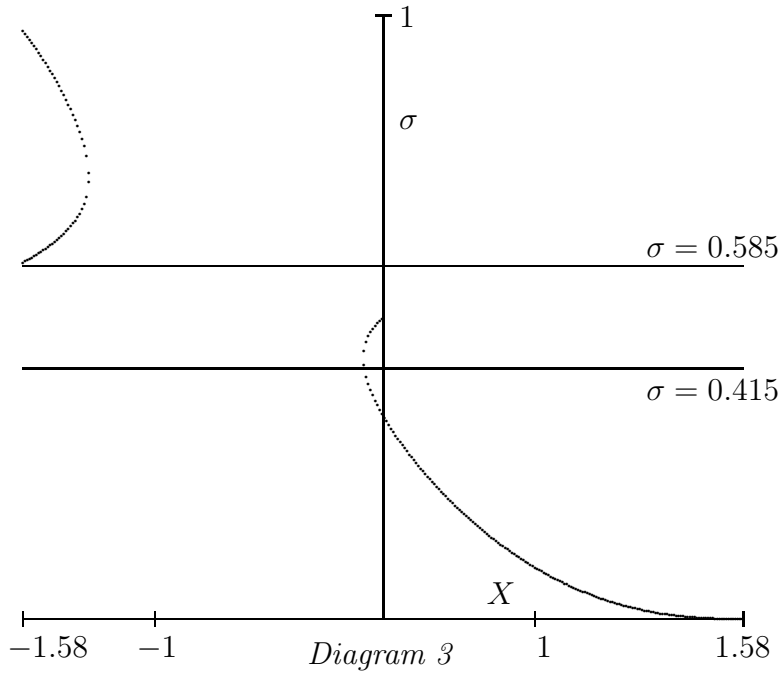
\input{Archpict_3}

%\unitlength1cm
%\begin{picture}(12,10) %Breite und Höhe

% Koordinaten wie üblich: (11,1) heißt im Bild 11(cm) nach rechts und 1(cm)
% nach oben, bei Geraden ist (a,b) der Anstieg b/a und {c} die Länge,

%\put(1,1){\line(1,0){10}}
%\put(1,1){\line(0,1){8.5}}
%\put(11,1){\line(0,1){8.5}}
%\put(1,9.5){\line(1,0){10}}
%\put(4,0.5){\hspace{1cm} {\em Diagram 3}}

%\end{picture}

The considerable differences between our examples suggest that a {\em global} theory of the
non-archimedean case (such as a theorem about the number of equilibria for a given pair
$(a,\sigma)$) may be a difficult matter.

%Within these bounds, however, an archimedean solution is possible only for $b=0$,
%because of (\ref{1.12}).

%%%%%%%%%%%%%%%%%%%%%%%%%%%%%%%%%%%%%%%%%%%%%%%%%%%%%%%%%%%%%%%%%
\section*{4. Classification of equilibria in the archimedean case}
%%%%%%%%%%%%%%%%%%%%%%%%%%%%%%%%%%%%%%%%%%%%%%%%%%%%%%%%%%%%%%%%%

The classification of the above equilibria requires considering
the potential energy of a certain position of the paraboloid segment $\PP$.
As in Section 1, we start with a right hand position in the archimedean case. Recall
that the center of gravity of the submerged part $\PP'$ is $B'=(x',0,z')$ with $x'=b/2$ and
$z'=5b^2/12+2c/3$. We need the moments
$x'V'$ and  $z'V'$, where
 $V'=a'^2\pi/2$ is the volume of $\PP'$ (with $a'$ as in (\ref{1.4})).
In what follows we use the abbreviation
\BE
\label{4.2}
  \beta=\sqrt{b^2+1}.
\EDE
Further, we work with a Hesse normal form of the plane
$\EE$ in order to describe the {\em distance} of a point $(x,y,z)$ from $\EE$,
namely
\BD
  \EE=\{(z-bx-c)/{\beta}=0\}.
\ED
Therefore, $-(z'-bx'-c)V'/\beta$ can be considered as the potential of the {\em buoyancy} of $\PP$.
In the same way the potential of the {\em weight} of $\PP$ is given by
$(2a/3-c)\sigma V/\beta$ since $\PP$ has the volume $V$ and the center of gravity
$B=(0,0,2a/3)$, see (\ref{1.2}).
Accordingly, our potential function
has the shape
\BD
U=U(c,b)=(({2a}/3-c)\sigma V+a'V'/3)/\beta
\ED
because $-z'+bx'+c=a'/3$. In this section we write
\BD
   F_0=V'-\sigma V \ENS \MB{ and }\ENS E_0=fV',
\ED
where $f=5b^2/12+2(c-a)/3+1/2$, see (\ref{1.7}).
The floating and equilibrium conditions read, thus, $F_0=0$ and $bE_0=0$,
respectively. One easily verifies
\BE
\label{4.3}
  \frac{\partial U}{\partial c\,}=\frac {F_0}{\beta},\ENS
  \frac{\partial U}{\partial b\,}=\frac b{\beta^3}\left(\left(\frac{2a}3-c\right)F_0+E_0\right).
\EDE
Equilibria should be the same as {\em stationary points} of the potential function. An
inspection of the derivatives of $U$ as given in (\ref{4.3}) proves this.
The second derivatives of $U$ can be written
\BE
\label{4.4}
  \frac{\partial^2 U}{\partial\, c^2}=\frac {2V'}{a'\beta},\ENS
  \frac{\partial^2 U}{\partial c\, \partial b}=\frac b{\beta}\left(\frac{V'}{a'}-\frac {F_0}{\beta^2}\right),
\EDE
and
\BE
\label{4.5}
\frac{\partial^2 U}{\partial\, b^2}=
\frac {b^2V'}{a'\beta^3}\left(\frac{5b^2}8+\frac{c+1}2\right)+\frac{1-2b^2}{\beta^5}
  \left(\left(\frac{2a}3-c\right)F_0+E_0\right).
\EDE
In the equilibrium case with $b=0$ we also have $F_0=0$. Therefore, (\ref{4.4}) and (\ref{4.5}) show
that the {\em Hesse matrix} of $U$ has the shape
\BD
   \left(
     \begin{array}{cc}
       2V'/a' & 0 \\
       0 & E_0 \\
     \end{array}
   \right).
\ED
Recalling $a'=a\sqrt{\sigma}>0$ we see that this matrix
is positive definite if, and only if, $E_0>0$, which is the same as saying $a<3/(4(1-\sqrt{\sigma}))$.
So in this case the stationary point of $U$ is a minimum and our right hand
equilibrium is {\em stable}; it becomes {\em unstable} (more precisely, a saddle point) for
$a>3/(4(1-\sqrt{\sigma}))$, whereas the case $a=3/(4(1-\sqrt{\sigma}))$ cannot be classified
in this way. Using higher derivatives one can show that this equilibrium is also stable
(as Archimedes did without this device, see \cite{Ar}, Lib. II, Sect. iv).
Similarly, the case $b<0$, $F_0=E_0=0$, gives the Hesse matrix
\BD
   \frac{V'}{a'\beta}
   \left(
     \begin{array}{cc}
       2 & b \\
       b & \left({5b^4}/8+(c+1)\,b^2/2\right)/{\beta^2}\\
     \end{array}
   \right).
\ED
Again, $a'=a\sqrt{\sigma}>0$. Since the determinant of this matrix equals
$V'^2b^2/(a'\beta^4)>0$, we see that equilibria with $b\ne0$ are {\em stable}.

\MN{\em Remarks.} 1. In our model no {\em maximum} of the potential function is to be expected,
since the potential energy will always grow if one moves  $\PP$ \BQ upwards\EQ,
i.e., more to the right hand side and, simultaneously, in a direction perpendicular to $\EE$.
This observation applies to left hand positions in an analogous way, for instance, to the non-archimedean
positions of the next section.

2. The {\em left hand} position that corresponds to our right hand one has the same potential,
see \cite{Gi} and the end of Section 1. This can be checked directly if one observes
that $\sigma^*=1-\sigma$ plays the role of $\sigma$ and
\BD
  B''=B+\frac{V'}{V-V'}(B-B')
\ED
that of $B$; observe, further, that the weight and buoyancy potentials change their signs.

%%%%%%%%%%%%%%%%%%%%%%%%%%%%%%%%%%%%%%%%%%%%%%%%%%%%%%%%%%%%%%%%%
\section*{5. Classification of equilibria in the non-archimedean case}
%%%%%%%%%%%%%%%%%%%%%%%%%%%%%%%%%%%%%%%%%%%%%%%%%%%%%%%%%%%%%%%%%

As in Section 2, we consider a left hand position of $\PP$ in the non-archimedean case
and adopt the corresponding notations.
In order to define the potential function $U$, we use the same Hesse normal form of $\EE$ as
as in the foregoing section. Then the potential of the weight of $\PP$ remains unchanged up to the
sign, i.e., it equals
$-(2a/3-c)\sigma V/\beta$, with $V$ as in (\ref{1.2}) and $\beta$ as in (\ref{4.2}).
The sign change is due to the transition from
a right hand position to a left hand one. Similarly, the potential of the buoyancy is
$(z'-bx'-c)(V_1-V_2)/\beta$.
Hence we have
\BD
   U=\frac 1{\beta}\left(\left(z'-bx'-c\right)(V_1-V_2)-\left(\frac{2a}3-c\right)\sigma V\right).
\ED
In the spirit of Section 2, we consider $U$
as a function of $X$ and $b$ instead of $c$ and $b$, so we write $c=a-bX$. By means of
formulas (\ref{2.13}) and (\ref{2.12}) we obtain
\BE
\label{5.2}
 U=U(X,b)= \frac 1{\beta}\left(\left(\frac a3-bX\right)
  \left(\sigma V-V_1\right)+\frac{a'}3 V_2-\frac{2bA^{5/2}}9\right)
\EDE
(recall (\ref{2.4})). Now the analogue of (\ref{4.3}) reads
\BE
\label{5.3}
  \frac{\partial U}{\partial X\,}=\frac {bF}{\beta},\ENS
  \frac{\partial U}{\partial b\,}=\frac 1{\beta^3}\left(bE+\left(X+\frac{ba}3\right)F\right)
\EDE
with $E$ and $F$ as in Theorem \ref{t1}.
We also note the analogues of (\ref{4.4}) and (\ref{4.5}), namely,
\BE
\label{5.4}
  \frac{\partial^2 U}{\partial\, X^2}=\frac {2b^2}{3a'\beta}(3V_2+X'A^{3/2}),\ENS
  \frac{\partial^2 U}{\partial X\, \partial b}=
  \frac {b^2+6}{4a'\beta}V_2+\frac{F}{\beta^3}-\frac {3E}{a'\beta},
\EDE
where $X'=X-b/2$ (see (\ref{2.4.1})) and
\BE
\label{5.5}
\frac{\partial^2 U}{\partial\, b^2}=
\frac 1{8a'b\beta^3}(-2Xb^4+(4a-7)b^3+14Xb^2-6b+12X)V_2
+F_1+E_1,
\EDE
with
\BEA
\label{5.6}
F_1&=&\frac 1{3\beta^5} (-2ab^2-9Xb+a)\,F\ENS \MB{ and }\nonumber\\
E_1&=& \frac 1{4a'b\beta^5}(4b^5-4Xb^4+(13-8a)b^3-28Xb^2+(4a+6)b-12X)\,E\,.
\EEA
Of course, the correctness of these formulas is easy to check by the aid of a computer algebra system,
say. However, it takes some effort to find them and we think, therefore, that it is justified
to render them here.
In the case of an equilibrium we have $F_1=E_1=0$, so (\ref{5.6}) disappears and (\ref{5.5}) looks
fairly simple then. This obviously happens for the second item of (\ref{5.4}), too.
Further, formulas (\ref{5.3}) show that there is hardly a simpler
characterization of stationary points of $U$
than our floating and equilibrium conditions.

\MN
{\em Example 3.} We return to Example 1. Our formulas (\ref{5.4}) and (\ref{5.5})
quickly give the Hesse matrix of the equilibrium positions we described there.
From altogether three left hand equilibria (with $\sigma=0.51$)
the first and the third one (with tilt angles of about $131.653^{\circ}$ and
$94.506^{\circ}$) are stable, since the Hesse matrices have the pairs of approximate eigenvalues
$(0.00101514$, $6.83907084)$ and $(0.00001567$, $7.50021176)$ in these cases. The second left hand
equilibrium ($131.260^{\circ}$) is unstable (a saddle point), the respective eigenvalues being approximately
$(-0.00098808, 6.78938522)$. In the same way the first of the right hand equilibria ($34.961^{\circ}$)
is stable and the second one ($56.793^{\circ}$) a saddle point. The case when the first and the second
equilibrium melt into one (with an angle of about $131.456^{\circ}$)
can be settled by the aid of higher
derivatives. For this purpose we write
$(X_0,b_0)$ for the corresponding value $\approx (-1.02702703,-1.13205421)$ of $(X,b)$.
We use the substitution $X=Y+\lambda b$, where
\BD
 \lambda=-\frac{{\partial^2 U}/{\partial X\partial b}}{{\partial^2 U}/{\partial\,X^2}}(X_0,b_0).
\ED
On differentiating $U$ with respect to $Y$ and $b$, one sees that the second derivatives vanish for the
respective point with the exception of
${\partial^2 U}/{\partial\,Y^2}$, whereas ${\partial^3 U}/{\partial\,b^3}$ takes a value
$c \approx 0.20378903$. Hence the behaviour of $U$ in a neighbourhood of $(X_0,b_0)$
is like that of $cZ^3+O(Y^2+|Y|Z^2+Z^4)$ for $Z,Y$ close to zero.
This expression, however, becomes negative for $Y=0$
and $Z<0$, $|Z|$ small. Therefore, $(X_0,b_0)$ defines an unstable position.

Of course, we can follow the same line when we investigate the
stability properties of those points on the curve $(X,\sigma)$
that are rendered in Diagram 2. We find that all points on the upper branch describe stable
equilibria. As to the lower branch, those with $X\le -1.03$ belong to stable equilibria
and the remaining ones (with $X\ge -1.02$) to saddle points. It seems, thus, that
$(X_0,\sigma_0)$ with $\sigma_0\approx 0.51000554$ (which comes from the above pair $(X_0,b_0)$)
forms a sort of limit point for the stability of left hand equilibria.
Additional computations (with smaller steps) confirm this observation.
%as does its supplement of $\approx 48.544^{\circ}$
%in the case of right hand ones.

Again, the behaviour of Example 2 is different. All points of Diagram 3 with $X<-1$ (i.e.,
those of the left component of the curve) give stable equilibria, whereas the remaining
ones belong to saddle points.

%%%%%%%%%%%%%%%%%%%%%%%%%%%%%%%%%%%%%%%%%%%%%%%%%%%%%%%%%%%%%%%%%
\section*{6. The horizontal case}
%%%%%%%%%%%%%%%%%%%%%%%%%%%%%%%%%%%%%%%%%%%%%%%%%%%%%%%%%%%%%%%%%

The case when the axis of the paraboloid segment is horizontal was not treated so far.
It will henceforth be called the {\em horizontal case}.
Our attempts to include this case in the above discussion gave rise to problems with differentiability,
so it seems justified not to do so.
In our setting the horizontal case can be characterized by  $\EE=\{x=X\}$,
where $X$ satisfies $-\sqrt a<X<\sqrt a$, see (\ref{2.1}).
One easily checks that the equilibrium condition takes the simple shape
 $(2a/3-z_1)V_1=0$. By (\ref{2.9}), this is equivalent
to $4XA^{5/2}/45=0$ and, therefore, to $X=0$. In this case $V_1=V/2$ holds for the
respective volumes, so the floating condition $V_1=\sigma V$ requires $\sigma=1/2$.
Accordingly, only this rather obvious equilibrium position is possible in the horizontal
case. The {\em classification} of this equilibrium, however, is less obvious.

To this end we use the potential $U_0$ of this position; it takes the
value $4a ^{5/2}/15$, as is readily seen.
Moreover, we consider
neighbouring {\em left hand} positions, i. e., pairs $(X,b)$ with $|X|$ small and $-b$
large. It is advantageous to work with $X'=X-b/2\gg 0$ instead of $b$ (see (\ref{2.4.1})).
Therefore, the potential function $U$ of (\ref{5.2}) reads $U=U(X,X')$ now.
Then we insert the series
\BD
 \arctan \frac{X}{\sqrt A}=\frac X{\sqrt A}-\frac{X^3}{3A^{3/2}}+\frac{X^5}{5A^{5/2}}-\LD,\ENS
 \arctan \frac{X'}{\sqrt A}=\frac{\pi}2-\frac{\sqrt A}{X'}+ \frac{A^{3/2}}{3X'^3}-\LD,
\ED
and
\BD
 \frac{1}{\beta}=\frac{1}{\sqrt{4(X'-X)^2+1}}=\frac 1{2X'}+\frac {X}{2X'^2}+ \frac{X^2/2-1/16}{X'^3}-\LD
\ED
together with
\BD
\frac{1}{A^{5/2}}=\frac 1{a^{5/2}}+\frac {5X^2}{2a^{7/2}}+ \frac{35X^4}{8a^{9/2}}-\LD,
\ED
into (\ref{5.2}) and obtain, in a straightforward (though laborious) way,
\BD
%\label{6.2}
U(X,X')=U_0+\frac{2a^{3/2}}3X^2+\frac{4a^{5/2}}{15}XY+
\left(\frac{4 a^{7/2}}{105}-\frac{a^{5/2}}{30}\right)Y^2+O((|X|+Y)^3),
\ED
where $Y=1/X'$ is positive and close to $0$.
The quadratic form in $X$, $Y$ on the right hand side is positive definite for $a>35/12$
and indefinite for $a<35/12$, thus indicating stability and instability, respectively.
To settle the case $a=35/12$, we need more terms of
this expansion and the substitution $Y=Y'-12X/7$. This gives
\BD
U=U_0-\frac{79\sqrt{105}}{3969}X^4+O(X^2|Y'|+Y'^2+|X|^5)
\ED
(where $U_0$ takes the value ${245\sqrt{105}}/{648}$).
So we can choose, for each small value of $|X|$, a
number $Y'$ with $|Y'|<|X|$ such that $U<U_0$. Hence this equilibrium is unstable.

%%%%%%%%%%%%%%%%%%%%%%%%%%%%%%%%%%%%%%%%%%%%%%%%%%%%%
%%%%%%%%%%%%%%%%%%%%%%%%%%%%%%%%%%%%%%%%%%%%%%%%%%%%%%%%%%%%%%%%%%%%%%%%%%

\vspace{0.5cm}
\noindent
Kurt Girstmair and
Gerhard Kirchner\\
Institut f\"ur Mathematik\\
Universit\"at Innsbruck\\
Technikerstr. 13/7\\
A-6020 Innsbruck, Austria\\
Kurt.Girstmair@uibk.ac.at,
Gehard.Kirchner@uibk.ac.at

\end{document}

%% file: Archpict_2.tex
%Varianten:\circle* und/oder {0.04}
\def\PKT{\circle*{0.01}}
\def\PT#1#2{\put(#1,#2){\PKT}}

%%%% Hier beginnt das äußere BILD!

\unitlength1cm
\begin{picture}(13,9.5) %Breite und Höhe

% Koordinaten wie üblich: (11,1) heißt im Bild 11(cm) nach rechts und 1(cm)
% nach oben, bei Geraden ist (a,b) der Anstieg b/a und {c} die Länge,

%\put(1,1){\line(1,0){10}}
%\put(1,1){\line(0,1){8.5}}
%\put(11,1){\line(0,1){8.5}}
%\put(1,9.5){\line(1,0){10}}
\put(4.3,0.3){\hspace{1cm} {\em Diagram 2}}

%Inneres Bild
%***************************************
\put(1.2,0.9){\begin{picture}(12,8)
%***************************************
\put(0,0){\line(1,0){10.68}}
\put(0,4.08){\line(1,0){10.68}}
\put(0,3.92){\line(1,0){10.68}}
\put(5.34,0){\line(0,1){8}}
\put(0,-0.1){\line(0,1){0.2}}
\put(-0.3,-0.5){$-1.78$}
\put(2.34,-0.1){\line(0,1){0.2}}
\put(2.1,-0.5){$-1$}
\put(8.34,-0.1){\line(0,1){0.2}}
\put(8.3,-0.5){$1$}
\put(10.68,-0.1){\line(0,1){0.2}}
\put(10.3,-0.5){$1.78$}
\put(6.5,0.2){$X$}
\put(9.6,3.5){$\sigma=0.49$}
\put(9.6,4.25){$\sigma=0.51$}
\put(5.24,8){\line(1,0){0.2}}
\put(5.6,7.85){$1$}
\put(5.6,6.5){$\sigma$}

\PT{  0.000}{ 7.883}    % 35  0.00000000000009
\PT{  0.030}{ 7.860}    % 31 -0.00000000000060
\PT{  0.060}{ 7.837}    % 31 -0.00000000000046
\PT{  0.090}{ 7.813}    % 31 -0.00000000000026
\PT{  0.120}{ 7.788}    % 32  0.00000000000015
\PT{  0.150}{ 7.763}    % 32 -0.00000000000024
\PT{  0.180}{ 7.737}    % 29  0.00000000000030
\PT{  0.210}{ 7.711}    % 32  0.00000000000071
\PT{  0.240}{ 7.684}    % 31 -0.00000000000015
\PT{  0.270}{ 7.657}    % 32  0.00000000000005
\PT{  0.300}{ 7.629}    % 30  0.00000000000063
\PT{  0.330}{ 7.601}    % 30  0.00000000000034
\PT{  0.360}{ 7.573}    % 29  0.00000000000100
\PT{  0.390}{ 7.544}    % 31 -0.00000000000005
\PT{  0.420}{ 7.515}    % 31 -0.00000000000092
\PT{  0.450}{ 7.486}    % 32 -0.00000000000005
\PT{  0.480}{ 7.456}    % 30 -0.00000000000003
\PT{  0.510}{ 7.426}    % 31 -0.00000000000032
\PT{  0.540}{ 7.396}    % 30  0.00000000000077
\PT{  0.570}{ 7.365}    % 27  0.00000000000051
\PT{  0.600}{ 7.335}    % 32 -0.00000000000019
\PT{  0.630}{ 7.304}    % 31 -0.00000000000059
\PT{  0.660}{ 7.272}    % 30  0.00000000000066
\PT{  0.690}{ 7.241}    % 28  0.00000000000027
\PT{  0.720}{ 7.209}    % 32 -0.00000000000038
\PT{  0.750}{ 7.177}    % 30 -0.00000000000037
\PT{  0.780}{ 7.145}    % 31  0.00000000000073
\PT{  0.810}{ 7.113}    % 30  0.00000000000066
\PT{  0.840}{ 7.081}    % 32 -0.00000000000002
\PT{  0.870}{ 7.048}    % 31 -0.00000000000094
\PT{  0.900}{ 7.015}    % 30  0.00000000000025
\PT{  0.930}{ 6.983}    % 30  0.00000000000025
\PT{  0.960}{ 6.949}    % 31  0.00000000000083
\PT{  0.990}{ 6.916}    % 31  0.00000000000067
\PT{  1.020}{ 6.883}    % 31 -0.00000000000048
\PT{  1.050}{ 6.850}    % 31 -0.00000000000071
\PT{  1.080}{ 6.816}    % 32  0.00000000000046
\PT{  1.110}{ 6.782}    % 30 -0.00000000000057
\PT{  1.140}{ 6.749}    % 31  0.00000000000065
\PT{  1.170}{ 6.715}    % 29  0.00000000000063
\PT{  1.200}{ 6.681}    % 29 -0.00000000000050
\PT{  1.230}{ 6.647}    % 31  0.00000000000003
\PT{  1.260}{ 6.612}    % 32  0.00000000000027
\PT{  1.290}{ 6.578}    % 32 -0.00000000000017
\PT{  1.320}{ 6.544}    % 31  0.00000000000016
\PT{  1.350}{ 6.510}    % 32  0.00000000000041
\PT{  1.380}{ 6.475}    % 31  0.00000000000031
\PT{  1.410}{ 6.441}    % 31  0.00000000000044
\PT{  1.440}{ 6.406}    % 28  0.00000000000097
\PT{  1.470}{ 6.372}    % 25  0.00000000000027
\PT{  1.500}{ 6.337}    % 27 -0.00000000000039
\PT{  1.530}{ 6.302}    % 31  0.00000000000066
\PT{  1.560}{ 6.268}    % 28 -0.00000000000029
\PT{  1.590}{ 6.233}    % 30  0.00000000000039
\PT{  1.620}{ 6.198}    % 28 -0.00000000000088
\PT{  1.650}{ 6.164}    % 29 -0.00000000000037
\PT{  1.680}{ 6.129}    % 31  0.00000000000069
\PT{  1.710}{ 6.094}    % 30 -0.00000000000024
\PT{  1.740}{ 6.059}    % 32 -0.00000000000030
\PT{  1.770}{ 6.025}    % 31 -0.00000000000090
\PT{  1.800}{ 5.990}    % 32  0.00000000000025
\PT{  1.830}{ 5.956}    % 32 -0.00000000000009
\PT{  1.860}{ 5.921}    % 29 -0.00000000000055
\PT{  1.890}{ 5.886}    % 31  0.00000000000090
\PT{  1.920}{ 5.852}    % 31 -0.00000000000037
\PT{  1.950}{ 5.817}    % 31  0.00000000000064
\PT{  1.980}{ 5.783}    % 30 -0.00000000000006
\PT{  2.010}{ 5.749}    % 31  0.00000000000012
\PT{  2.040}{ 5.714}    % 28 -0.00000000000048
\PT{  2.070}{ 5.680}    % 31  0.00000000000010
\PT{  2.100}{ 5.646}    % 28  0.00000000000008
\PT{  2.130}{ 5.612}    % 31  0.00000000000033
\PT{  2.160}{ 5.578}    % 31 -0.00000000000051
\PT{  2.190}{ 5.545}    % 31  0.00000000000091
\PT{  2.220}{ 5.511}    % 28  0.00000000000038
\PT{  2.250}{ 5.478}    % 28 -0.00000000000071
\PT{  2.280}{ 5.444}    % 31  0.00000000000051
\PT{  2.310}{ 5.411}    % 28 -0.00000000000023
\PT{  2.340}{ 5.378}    % 31 -0.00000000000016
\PT{  2.370}{ 5.345}    % 31  0.00000000000049
\PT{  2.400}{ 5.313}    % 30  0.00000000000036
\PT{  2.430}{ 5.280}    % 30 -0.00000000000098
\PT{  2.460}{ 5.248}    % 30  0.00000000000049
\PT{  2.490}{ 5.216}    % 30 -0.00000000000092
\PT{  2.520}{ 5.185}    % 29 -0.00000000000091
\PT{  2.550}{ 5.153}    % 30  0.00000000000020
\PT{  2.580}{ 5.122}    % 25  0.00000000000054
\PT{  2.610}{ 5.091}    % 30 -0.00000000000077
\PT{  2.640}{ 5.061}    % 29 -0.00000000000038
\PT{  2.670}{ 5.031}    % 31  0.00000000000004
\PT{  2.700}{ 5.001}    % 30  0.00000000000063
\PT{  2.730}{ 4.972}    % 30  0.00000000000003
\PT{  2.760}{ 4.943}    % 30  0.00000000000054
\PT{  2.790}{ 4.915}    % 29  0.00000000000037
\PT{  2.820}{ 4.887}    % 30  0.00000000000071
\PT{  2.850}{ 4.860}    % 29 -0.00000000000030
\PT{  2.880}{ 4.834}    % 29  0.00000000000033
\PT{  2.910}{ 4.808}    % 30 -0.00000000000077
\PT{  2.940}{ 4.782}    % 30 -0.00000000000044
\PT{  2.970}{ 4.757}    % 27  0.00000000000062
\PT{  3.000}{ 4.733}    % 31 -0.00000000000074
\PT{  3.030}{ 4.710}    % 30  0.00000000000040
\PT{  3.060}{ 4.687}    % 31 -0.00000000000065
\PT{  3.090}{ 4.665}    % 30  0.00000000000025
\PT{  3.120}{ 4.644}    % 31  0.00000000000020
\PT{  3.150}{ 4.624}    % 31 -0.00000000000068
\PT{  3.180}{ 4.605}    % 29  0.00000000000023
\PT{  3.210}{ 4.586}    % 30 -0.00000000000056
\PT{  3.240}{ 4.568}    % 28  0.00000000000005
\PT{  3.270}{ 4.551}    % 30  0.00000000000040
\PT{  3.300}{ 4.534}    % 29 -0.00000000000085
\PT{  3.330}{ 4.519}    % 32  0.00000000000005
\PT{  3.360}{ 4.504}    % 29  0.00000000000021
\PT{  3.390}{ 4.489}    % 31  0.00000000000068
\PT{  3.420}{ 4.476}    % 30 -0.00000000000064
\PT{  3.450}{ 4.463}    % 26 -0.00000000000087
\PT{  3.480}{ 4.450}    % 28  0.00000000000061
\PT{  3.510}{ 4.438}    % 31  0.00000000000048
\PT{  3.540}{ 4.427}    % 32 -0.00000000000021
\PT{  3.570}{ 4.416}    % 32 -0.00000000000036
\PT{  3.600}{ 4.406}    % 30  0.00000000000015
\PT{  3.630}{ 4.396}    % 28 -0.00000000000017
\PT{  3.660}{ 4.386}    % 31 -0.00000000000039
\PT{  3.690}{ 4.376}    % 29  0.00000000000093
\PT{  3.720}{ 4.367}    % 31  0.00000000000048
\PT{  3.750}{ 4.359}    % 28  0.00000000000035
\PT{  3.780}{ 4.350}    % 32 -0.00000000000062
\PT{  3.810}{ 4.342}    % 31 -0.00000000000057
\PT{  3.840}{ 4.334}    % 32 -0.00000000000080
\PT{  3.870}{ 4.326}    % 32 -0.00000000000061
\PT{  3.900}{ 4.319}    % 29  0.00000000000023
\PT{  3.930}{ 4.311}    % 32  0.00000000000099
\PT{  3.960}{ 4.304}    % 30 -0.00000000000031
\PT{  3.990}{ 4.297}    % 32  0.00000000000068
\PT{  4.020}{ 4.289}    % 32  0.00000000000015
\PT{  4.050}{ 4.282}    % 33  0.00000000000013
\PT{  4.080}{ 4.275}    % 32  0.00000000000021
\PT{  4.110}{ 4.269}    % 31  0.00000000000025
\PT{  4.140}{ 4.262}    % 33  0.00000000000024
\PT{  4.170}{ 4.255}    % 31  0.00000000000036
\PT{  4.200}{ 4.248}    % 33 -0.00000000000027
\PT{  4.230}{ 4.242}    % 29 -0.00000000000029
\PT{  4.260}{ 4.235}    % 29  0.00000000000056
\PT{  4.290}{ 4.229}    % 32  0.00000000000067
\PT{  4.320}{ 4.222}    % 32 -0.00000000000046
\PT{  4.350}{ 4.216}    % 32 -0.00000000000008
\PT{  4.380}{ 4.209}    % 31 -0.00000000000070
\PT{  4.410}{ 4.203}    % 33 -0.00000000000001
\PT{  4.440}{ 4.196}    % 31  0.00000000000061
\PT{  4.470}{ 4.190}    % 32 -0.00000000000086
\PT{  4.500}{ 4.183}    % 32 -0.00000000000095
\PT{  4.530}{ 4.177}    % 31 -0.00000000000092
\PT{  4.560}{ 4.170}    % 33 -0.00000000000054
\PT{  4.590}{ 4.164}    % 32  0.00000000000034
\PT{  4.620}{ 4.157}    % 32  0.00000000000046
\PT{  4.650}{ 4.151}    % 31 -0.00000000000049
\PT{  4.680}{ 4.145}    % 32  0.00000000000018
\PT{  4.710}{ 4.138}    % 33 -0.00000000000039
\PT{  4.740}{ 4.132}    % 31  0.00000000000083
\PT{  4.770}{ 4.125}    % 33 -0.00000000000021
\PT{  4.800}{ 4.119}    % 33 -0.00000000000001
\PT{  4.830}{ 4.112}    % 33  0.00000000000034
\PT{  4.860}{ 4.106}    % 33  0.00000000000007
\PT{  4.890}{ 4.099}    % 32  0.00000000000093
\PT{  4.920}{ 4.092}    % 33  0.00000000000002
\PT{  4.950}{ 4.086}    % 31  0.00000000000030
\PT{  4.980}{ 4.079}    % 33 -0.00000000000006
\PT{  5.010}{ 4.073}    % 33  0.00000000000066
\PT{  5.040}{ 4.066}    % 33  0.00000000000004
\PT{  5.070}{ 4.060}    % 32 -0.00000000000074
\PT{  5.100}{ 4.053}    % 33 -0.00000000000015
\PT{  5.130}{ 4.046}    % 32  0.00000000000094
\PT{  5.160}{ 4.040}    % 32  0.00000000000028
\PT{  5.190}{ 4.033}    % 31 -0.00000000000043
\PT{  5.220}{ 4.027}    % 29 -0.00000000000010
\PT{  5.250}{ 4.020}    % 31  0.00000000000046
\PT{  5.280}{ 4.013}    % 32  0.00000000000093
\PT{  5.310}{ 4.007}    % 30 -0.00000000000065

\PT{  0.000}{ 3.705}    % 40 -0.00000000000047
\PT{  0.030}{ 3.712}    % 34  0.00000000000011
\PT{  0.060}{ 3.718}    % 34 -0.00000000000053
\PT{  0.090}{ 3.725}    % 35 -0.00000000000057
\PT{  0.120}{ 3.732}    % 35  0.00000000000013
\PT{  0.150}{ 3.739}    % 34  0.00000000000073
\PT{  0.180}{ 3.745}    % 33  0.00000000000041
\PT{  0.210}{ 3.752}    % 34  0.00000000000075
\PT{  0.240}{ 3.759}    % 31  0.00000000000092
\PT{  0.270}{ 3.766}    % 35 -0.00000000000033
\PT{  0.300}{ 3.773}    % 35 -0.00000000000028
\PT{  0.330}{ 3.780}    % 33  0.00000000000055
\PT{  0.360}{ 3.788}    % 33 -0.00000000000052
\PT{  0.390}{ 3.795}    % 34  0.00000000000070
\PT{  0.420}{ 3.802}    % 31 -0.00000000000046
\PT{  0.450}{ 3.809}    % 35  0.00000000000026
\PT{  0.480}{ 3.816}    % 33 -0.00000000000061
\PT{  0.510}{ 3.823}    % 34  0.00000000000085
\PT{  0.540}{ 3.830}    % 35 -0.00000000000003
\PT{  0.570}{ 3.837}    % 34 -0.00000000000040
\PT{  0.600}{ 3.843}    % 33 -0.00000000000089
\PT{  0.630}{ 3.850}    % 34 -0.00000000000051
\PT{  0.660}{ 3.857}    % 32  0.00000000000078
\PT{  0.690}{ 3.864}    % 34  0.00000000000046
\PT{  0.720}{ 3.871}    % 32  0.00000000000017
\PT{  0.750}{ 3.877}    % 35  0.00000000000016
\PT{  0.780}{ 3.884}    % 34 -0.00000000000055
\PT{  0.810}{ 3.890}    % 30 -0.00000000000054
\PT{  0.840}{ 3.897}    % 33 -0.00000000000088
\PT{  0.870}{ 3.903}    % 35 -0.00000000000001
\PT{  0.900}{ 3.910}    % 32 -0.00000000000083
\PT{  0.930}{ 3.916}    % 34 -0.00000000000002
\PT{  0.960}{ 3.922}    % 34 -0.00000000000083
\PT{  0.990}{ 3.928}    % 33  0.00000000000026
\PT{  1.020}{ 3.935}    % 34 -0.00000000000068
\PT{  1.050}{ 3.940}    % 34  0.00000000000066
\PT{  1.080}{ 3.946}    % 33 -0.00000000000097
\PT{  1.110}{ 3.952}    % 34 -0.00000000000043
\PT{  1.140}{ 3.958}    % 32 -0.00000000000079
\PT{  1.170}{ 3.964}    % 34 -0.00000000000060
\PT{  1.200}{ 3.969}    % 31 -0.00000000000036
\PT{  1.230}{ 3.974}    % 33 -0.00000000000025
\PT{  1.260}{ 3.980}    % 33  0.00000000000040
\PT{  1.290}{ 3.985}    % 28  0.00000000000029
\PT{  1.320}{ 3.990}    % 33 -0.00000000000006
\PT{  1.350}{ 3.995}    % 34  0.00000000000036
\PT{  1.380}{ 4.000}    % 32 -0.00000000000060
\PT{  1.410}{ 4.005}    % 33 -0.00000000000002
\PT{  1.440}{ 4.010}    % 33  0.00000000000016
\PT{  1.470}{ 4.014}    % 33 -0.00000000000034
\PT{  1.500}{ 4.019}    % 32 -0.00000000000090
\PT{  1.530}{ 4.023}    % 29 -0.00000000000055
\PT{  1.560}{ 4.027}    % 23  0.00000000000065
\PT{  1.590}{ 4.031}    % 33 -0.00000000000099
\PT{  1.620}{ 4.035}    % 32 -0.00000000000048
\PT{  1.650}{ 4.039}    % 34 -0.00000000000001
\PT{  1.680}{ 4.042}    % 34 -0.00000000000027
\PT{  1.710}{ 4.046}    % 34  0.00000000000068
\PT{  1.740}{ 4.049}    % 34 -0.00000000000047
\PT{  1.770}{ 4.052}    % 34 -0.00000000000036
\PT{  1.800}{ 4.056}    % 34  0.00000000000006
\PT{  1.830}{ 4.058}    % 33 -0.00000000000090
\PT{  1.860}{ 4.061}    % 32 -0.00000000000038
\PT{  1.890}{ 4.064}    % 31  0.00000000000009
\PT{  1.920}{ 4.066}    % 32  0.00000000000087
\PT{  1.950}{ 4.068}    % 33 -0.00000000000084
\PT{  1.980}{ 4.070}    % 32 -0.00000000000080
\PT{  2.010}{ 4.072}    % 33 -0.00000000000094
\PT{  2.040}{ 4.074}    % 32 -0.00000000000042
\PT{  2.070}{ 4.075}    % 34  0.00000000000031
\PT{  2.100}{ 4.077}    % 28 -0.00000000000031
\PT{  2.130}{ 4.078}    % 33 -0.00000000000040
\PT{  2.160}{ 4.079}    % 34  0.00000000000018
\PT{  2.190}{ 4.079}    % 33  0.00000000000037
\PT{  2.220}{ 4.080}    % 33  0.00000000000030
\PT{  2.250}{ 4.080}    % 33 -0.00000000000028
\PT{  2.280}{ 4.080}    % 33 -0.00000000000029
\PT{  2.310}{ 4.080}    % 33  0.00000000000003
\PT{  2.340}{ 4.079}    % 33 -0.00000000000065
\PT{  2.370}{ 4.078}    % 33  0.00000000000045
\PT{  2.400}{ 4.077}    % 33 -0.00000000000036
\PT{  2.430}{ 4.076}    % 32 -0.00000000000081
\PT{  2.460}{ 4.074}    % 28  0.00000000000058
\PT{  2.490}{ 4.072}    % 31 -0.00000000000073
\PT{  2.520}{ 4.069}    % 31 -0.00000000000086
\PT{  2.550}{ 4.067}    % 32 -0.00000000000094
\PT{  2.580}{ 4.064}    % 33  0.00000000000006
\PT{  2.610}{ 4.060}    % 32  0.00000000000034
\PT{  2.640}{ 4.056}    % 33 -0.00000000000016
\PT{  2.670}{ 4.052}    % 33  0.00000000000017
\PT{  2.700}{ 4.047}    % 32 -0.00000000000099
\PT{  2.730}{ 4.041}    % 33 -0.00000000000016
\PT{  2.760}{ 4.036}    % 32 -0.00000000000088
\PT{  2.790}{ 4.029}    % 32  0.00000000000045
\PT{  2.820}{ 4.022}    % 32  0.00000000000044
\PT{  2.850}{ 4.015}    % 32 -0.00000000000040
\PT{  2.880}{ 4.007}    % 31  0.00000000000042
\PT{  2.910}{ 3.998}    % 32 -0.00000000000064
\PT{  2.940}{ 3.988}    % 28 -0.00000000000048
\PT{  2.970}{ 3.978}    % 32 -0.00000000000047
\PT{  3.000}{ 3.967}    % 31  0.00000000000062
\PT{  3.030}{ 3.955}    % 32  0.00000000000028
\PT{  3.060}{ 3.943}    % 32 -0.00000000000060
\PT{  3.090}{ 3.930}    % 31  0.00000000000007
\PT{  3.120}{ 3.915}    % 32  0.00000000000006
\PT{  3.150}{ 3.901}    % 32 -0.00000000000042
\PT{  3.180}{ 3.885}    % 30 -0.00000000000064
\PT{  3.210}{ 3.868}    % 32  0.00000000000041
\PT{  3.240}{ 3.851}    % 32 -0.00000000000042
\PT{  3.270}{ 3.833}    % 31  0.00000000000098
\PT{  3.300}{ 3.814}    % 31 -0.00000000000006
\PT{  3.330}{ 3.794}    % 30 -0.00000000000002
\PT{  3.360}{ 3.774}    % 30  0.00000000000059
\PT{  3.390}{ 3.753}    % 30 -0.00000000000033
\PT{  3.420}{ 3.731}    % 29  0.00000000000046
\PT{  3.450}{ 3.708}    % 30  0.00000000000083
\PT{  3.480}{ 3.685}    % 31  0.00000000000030
\PT{  3.510}{ 3.662}    % 28  0.00000000000079
\PT{  3.540}{ 3.638}    % 23 -0.00000000000058
\PT{  3.570}{ 3.613}    % 29  0.00000000000056
\PT{  3.600}{ 3.588}    % 28 -0.00000000000012
\PT{  3.630}{ 3.563}    % 31 -0.00000000000040
\PT{  3.660}{ 3.537}    % 30 -0.00000000000067
\PT{  3.690}{ 3.511}    % 30 -0.00000000000046
\PT{  3.720}{ 3.484}    % 31 -0.00000000000021
\PT{  3.750}{ 3.457}    % 31 -0.00000000000000
\PT{  3.780}{ 3.430}    % 30 -0.00000000000088
\PT{  3.810}{ 3.403}    % 29  0.00000000000078
\PT{  3.840}{ 3.376}    % 31  0.00000000000025
\PT{  3.870}{ 3.348}    % 25  0.00000000000032
\PT{  3.900}{ 3.320}    % 31  0.00000000000023
\PT{  3.930}{ 3.292}    % 30  0.00000000000015
\PT{  3.960}{ 3.264}    % 27 -0.00000000000038
\PT{  3.990}{ 3.236}    % 29 -0.00000000000024
\PT{  4.020}{ 3.208}    % 31 -0.00000000000057
\PT{  4.050}{ 3.180}    % 31  0.00000000000038
\PT{  4.080}{ 3.152}    % 31  0.00000000000076
\PT{  4.110}{ 3.123}    % 31  0.00000000000002
\PT{  4.140}{ 3.095}    % 31 -0.00000000000050
\PT{  4.170}{ 3.066}    % 28  0.00000000000095
\PT{  4.200}{ 3.038}    % 31  0.00000000000094
\PT{  4.230}{ 3.009}    % 31  0.00000000000034
\PT{  4.260}{ 2.981}    % 31 -0.00000000000081
\PT{  4.290}{ 2.953}    % 31 -0.00000000000042
\PT{  4.320}{ 2.924}    % 31  0.00000000000074
\PT{  4.350}{ 2.896}    % 31 -0.00000000000014
\PT{  4.380}{ 2.867}    % 31 -0.00000000000012
\PT{  4.410}{ 2.839}    % 31  0.00000000000082
\PT{  4.440}{ 2.811}    % 31  0.00000000000006
\PT{  4.470}{ 2.783}    % 31 -0.00000000000039
\PT{  4.500}{ 2.755}    % 29  0.00000000000018
\PT{  4.530}{ 2.727}    % 30 -0.00000000000018
\PT{  4.560}{ 2.699}    % 30 -0.00000000000004
\PT{  4.590}{ 2.671}    % 28  0.00000000000016
\PT{  4.620}{ 2.643}    % 32  0.00000000000003
\PT{  4.650}{ 2.615}    % 31  0.00000000000001
\PT{  4.680}{ 2.587}    % 30 -0.00000000000061
\PT{  4.710}{ 2.560}    % 29  0.00000000000064
\PT{  4.740}{ 2.532}    % 30 -0.00000000000099
\PT{  4.770}{ 2.505}    % 32  0.00000000000020
\PT{  4.800}{ 2.478}    % 30  0.00000000000013
\PT{  4.830}{ 2.451}    % 31  0.00000000000047
\PT{  4.860}{ 2.424}    % 31 -0.00000000000085
\PT{  4.890}{ 2.397}    % 31  0.00000000000050
\PT{  4.920}{ 2.370}    % 28  0.00000000000001
\PT{  4.950}{ 2.343}    % 28 -0.00000000000002
\PT{  4.980}{ 2.316}    % 31 -0.00000000000080
\PT{  5.010}{ 2.290}    % 26  0.00000000000078
\PT{  5.040}{ 2.263}    % 32 -0.00000000000023
\PT{  5.070}{ 2.237}    % 31  0.00000000000090
\PT{  5.100}{ 2.211}    % 31 -0.00000000000060
\PT{  5.130}{ 2.185}    % 32 -0.00000000000026
\PT{  5.160}{ 2.159}    % 31  0.00000000000072
\PT{  5.190}{ 2.133}    % 28 -0.00000000000051
\PT{  5.220}{ 2.108}    % 28 -0.00000000000012
\PT{  5.250}{ 2.082}    % 31  0.00000000000095
\PT{  5.280}{ 2.057}    % 30 -0.00000000000062
\PT{  5.310}{ 2.032}    % 31  0.00000000000016
\PT{  5.340}{ 2.007}    % 31 -0.00000000000042
\PT{  5.370}{ 1.982}    % 28  0.00000000000024
\PT{  5.400}{ 1.957}    % 32  0.00000000000003
\PT{  5.430}{ 1.932}    % 31 -0.00000000000074
\PT{  5.460}{ 1.908}    % 30  0.00000000000100
\PT{  5.490}{ 1.883}    % 31  0.00000000000021
\PT{  5.520}{ 1.859}    % 31 -0.00000000000027
\PT{  5.550}{ 1.835}    % 31 -0.00000000000011
\PT{  5.580}{ 1.811}    % 31 -0.00000000000082
\PT{  5.610}{ 1.787}    % 26  0.00000000000003
\PT{  5.640}{ 1.763}    % 30 -0.00000000000089
\PT{  5.670}{ 1.740}    % 31 -0.00000000000013
\PT{  5.700}{ 1.716}    % 30 -0.00000000000062
\PT{  5.730}{ 1.693}    % 28  0.00000000000011
\PT{  5.760}{ 1.670}    % 31  0.00000000000093
\PT{  5.790}{ 1.647}    % 30 -0.00000000000099
\PT{  5.820}{ 1.624}    % 30 -0.00000000000082
\PT{  5.850}{ 1.602}    % 30  0.00000000000065
\PT{  5.880}{ 1.579}    % 30  0.00000000000014
\PT{  5.910}{ 1.557}    % 30 -0.00000000000022
\PT{  5.940}{ 1.535}    % 31 -0.00000000000024
\PT{  5.970}{ 1.513}    % 32 -0.00000000000011
\PT{  6.000}{ 1.491}    % 29 -0.00000000000073
\PT{  6.030}{ 1.469}    % 25 -0.00000000000069
\PT{  6.060}{ 1.448}    % 29 -0.00000000000089
\PT{  6.090}{ 1.426}    % 32  0.00000000000020
\PT{  6.120}{ 1.405}    % 31  0.00000000000044
\PT{  6.150}{ 1.384}    % 32  0.00000000000027
\PT{  6.180}{ 1.363}    % 27 -0.00000000000024
\PT{  6.210}{ 1.342}    % 30  0.00000000000038
\PT{  6.240}{ 1.322}    % 30 -0.00000000000090
\PT{  6.270}{ 1.301}    % 32  0.00000000000007
\PT{  6.300}{ 1.281}    % 29 -0.00000000000074
\PT{  6.330}{ 1.261}    % 31 -0.00000000000067
\PT{  6.360}{ 1.241}    % 32 -0.00000000000012
\PT{  6.390}{ 1.221}    % 28 -0.00000000000038
\PT{  6.420}{ 1.202}    % 27  0.00000000000039
\PT{  6.450}{ 1.182}    % 32 -0.00000000000022
\PT{  6.480}{ 1.163}    % 31 -0.00000000000092
\PT{  6.510}{ 1.144}    % 31 -0.00000000000079
\PT{  6.540}{ 1.125}    % 27  0.00000000000016
\PT{  6.570}{ 1.106}    % 31 -0.00000000000076
\PT{  6.600}{ 1.087}    % 30 -0.00000000000077
\PT{  6.630}{ 1.069}    % 32  0.00000000000006
\PT{  6.660}{ 1.051}    % 30  0.00000000000099
\PT{  6.690}{ 1.032}    % 31  0.00000000000014
\PT{  6.720}{ 1.015}    % 30 -0.00000000000001
\PT{  6.750}{ 0.997}    % 32 -0.00000000000013
\PT{  6.780}{ 0.979}    % 30 -0.00000000000033
\PT{  6.810}{ 0.962}    % 31  0.00000000000040
\PT{  6.840}{ 0.944}    % 31 -0.00000000000023
\PT{  6.870}{ 0.927}    % 31  0.00000000000000
\PT{  6.900}{ 0.910}    % 30 -0.00000000000044
\PT{  6.930}{ 0.893}    % 31 -0.00000000000051
\PT{  6.960}{ 0.877}    % 26  0.00000000000091
\PT{  6.990}{ 0.860}    % 31  0.00000000000051
\PT{  7.020}{ 0.844}    % 27 -0.00000000000046
\PT{  7.050}{ 0.828}    % 30 -0.00000000000006
\PT{  7.080}{ 0.812}    % 32 -0.00000000000001
\PT{  7.110}{ 0.796}    % 29 -0.00000000000076
\PT{  7.140}{ 0.780}    % 30  0.00000000000051
\PT{  7.170}{ 0.765}    % 30  0.00000000000035
\PT{  7.200}{ 0.750}    % 31 -0.00000000000010
\PT{  7.230}{ 0.735}    % 30 -0.00000000000007
\PT{  7.260}{ 0.720}    % 31 -0.00000000000001
\PT{  7.290}{ 0.705}    % 29  0.00000000000050
\PT{  7.320}{ 0.690}    % 29 -0.00000000000026
\PT{  7.350}{ 0.676}    % 31 -0.00000000000022
\PT{  7.380}{ 0.661}    % 28 -0.00000000000034
\PT{  7.410}{ 0.647}    % 31 -0.00000000000047
\PT{  7.440}{ 0.633}    % 31 -0.00000000000010
\PT{  7.470}{ 0.619}    % 27  0.00000000000035
\PT{  7.500}{ 0.606}    % 26 -0.00000000000071
\PT{  7.530}{ 0.592}    % 30  0.00000000000027
\PT{  7.560}{ 0.579}    % 31 -0.00000000000029
\PT{  7.590}{ 0.566}    % 30  0.00000000000080
\PT{  7.620}{ 0.553}    % 31  0.00000000000037
\PT{  7.650}{ 0.540}    % 29 -0.00000000000039
\PT{  7.680}{ 0.528}    % 28  0.00000000000039
\PT{  7.710}{ 0.515}    % 31 -0.00000000000029
\PT{  7.740}{ 0.503}    % 30  0.00000000000020
\PT{  7.770}{ 0.491}    % 31  0.00000000000048
\PT{  7.800}{ 0.479}    % 29 -0.00000000000011
\PT{  7.830}{ 0.467}    % 30 -0.00000000000063
\PT{  7.860}{ 0.455}    % 31  0.00000000000046
\PT{  7.890}{ 0.444}    % 31 -0.00000000000021
\PT{  7.920}{ 0.432}    % 30 -0.00000000000060
\PT{  7.950}{ 0.421}    % 30 -0.00000000000092
\PT{  7.980}{ 0.410}    % 31  0.00000000000014
\PT{  8.010}{ 0.399}    % 28 -0.00000000000088
\PT{  8.040}{ 0.389}    % 30  0.00000000000053
\PT{  8.070}{ 0.378}    % 28  0.00000000000068
\PT{  8.100}{ 0.368}    % 26 -0.00000000000016
\PT{  8.130}{ 0.358}    % 28  0.00000000000039
\PT{  8.160}{ 0.348}    % 31  0.00000000000010
\PT{  8.190}{ 0.338}    % 30  0.00000000000039
\PT{  8.220}{ 0.328}    % 27  0.00000000000061
\PT{  8.250}{ 0.319}    % 29 -0.00000000000032
\PT{  8.280}{ 0.309}    % 29  0.00000000000078
\PT{  8.310}{ 0.300}    % 30  0.00000000000004
\PT{  8.340}{ 0.291}    % 31  0.00000000000001
\PT{  8.370}{ 0.282}    % 30  0.00000000000045
\PT{  8.400}{ 0.273}    % 30 -0.00000000000064
\PT{  8.430}{ 0.265}    % 30 -0.00000000000025
\PT{  8.460}{ 0.256}    % 28  0.00000000000087
\PT{  8.490}{ 0.248}    % 30  0.00000000000002
\PT{  8.520}{ 0.240}    % 29  0.00000000000043
\PT{  8.550}{ 0.232}    % 30  0.00000000000007
\PT{  8.580}{ 0.224}    % 30 -0.00000000000060
\PT{  8.610}{ 0.217}    % 29 -0.00000000000051
\PT{  8.640}{ 0.209}    % 30  0.00000000000029
\PT{  8.670}{ 0.202}    % 29 -0.00000000000082
\PT{  8.700}{ 0.194}    % 30  0.00000000000053
\PT{  8.730}{ 0.187}    % 30 -0.00000000000042
\PT{  8.760}{ 0.180}    % 28 -0.00000000000045
\PT{  8.790}{ 0.174}    % 29 -0.00000000000092
\PT{  8.820}{ 0.167}    % 27  0.00000000000076
\PT{  8.850}{ 0.161}    % 30 -0.00000000000024
\PT{  8.880}{ 0.154}    % 30 -0.00000000000003
\PT{  8.910}{ 0.148}    % 29  0.00000000000001
\PT{  8.940}{ 0.142}    % 27 -0.00000000000034
\PT{  8.970}{ 0.136}    % 27 -0.00000000000065
\PT{  9.000}{ 0.130}    % 26  0.00000000000076
\PT{  9.030}{ 0.125}    % 28 -0.00000000000014
\PT{  9.060}{ 0.119}    % 29  0.00000000000082
\PT{  9.090}{ 0.114}    % 29 -0.00000000000055
\PT{  9.120}{ 0.109}    % 29  0.00000000000035
\PT{  9.150}{ 0.104}    % 28 -0.00000000000055
\PT{  9.180}{ 0.099}    % 29  0.00000000000012
\PT{  9.210}{ 0.094}    % 26  0.00000000000041
\PT{  9.240}{ 0.090}    % 27  0.00000000000095
\PT{  9.270}{ 0.085}    % 28  0.00000000000020
\PT{  9.300}{ 0.081}    % 27 -0.00000000000032
\PT{  9.330}{ 0.077}    % 28 -0.00000000000057
\PT{  9.360}{ 0.072}    % 29 -0.00000000000006
\PT{  9.390}{ 0.069}    % 29 -0.00000000000018
\PT{  9.420}{ 0.065}    % 27  0.00000000000004
\PT{  9.450}{ 0.061}    % 28 -0.00000000000086
\PT{  9.480}{ 0.057}    % 27 -0.00000000000066
\PT{  9.510}{ 0.054}    % 28 -0.00000000000025
\PT{  9.540}{ 0.051}    % 28  0.00000000000021
\PT{  9.570}{ 0.047}    % 28 -0.00000000000060
\PT{  9.600}{ 0.044}    % 28  0.00000000000052
\PT{  9.630}{ 0.041}    % 26  0.00000000000049
\PT{  9.660}{ 0.039}    % 26  0.00000000000085
\PT{  9.690}{ 0.036}    % 27  0.00000000000098
\PT{  9.720}{ 0.033}    % 28  0.00000000000010
\PT{  9.750}{ 0.031}    % 25 -0.00000000000010
\PT{  9.780}{ 0.028}    % 27 -0.00000000000072
\PT{  9.810}{ 0.026}    % 27 -0.00000000000048
\PT{  9.840}{ 0.024}    % 26  0.00000000000061
\PT{  9.870}{ 0.022}    % 26  0.00000000000085
\PT{  9.900}{ 0.020}    % 26 -0.00000000000080
\PT{  9.930}{ 0.018}    % 27 -0.00000000000036
\PT{  9.960}{ 0.016}    % 27 -0.00000000000033
\PT{  9.990}{ 0.015}    % 24 -0.00000000000098
\PT{ 10.020}{ 0.013}    % 27 -0.00000000000002
\PT{ 10.050}{ 0.012}    % 25  0.00000000000027
\PT{ 10.080}{ 0.011}    % 25  0.00000000000040
\PT{ 10.110}{ 0.009}    % 25  0.00000000000018
\PT{ 10.140}{ 0.008}    % 24 -0.00000000000028
\PT{ 10.170}{ 0.007}    % 26 -0.00000000000006
\PT{ 10.200}{ 0.006}    % 24  0.00000000000094
\PT{ 10.230}{ 0.005}    % 19  0.00000000000064
\PT{ 10.260}{ 0.004}    % 24 -0.00000000000067
\PT{ 10.290}{ 0.004}    % 24  0.00000000000035
\PT{ 10.320}{ 0.003}    % 24  0.00000000000063
\PT{ 10.350}{ 0.002}    % 24  0.00000000000041
\PT{ 10.380}{ 0.002}    % 24 -0.00000000000017
\PT{ 10.410}{ 0.001}    % 22  0.00000000000005
\PT{ 10.440}{ 0.001}    % 23  0.00000000000043
\PT{ 10.470}{ 0.001}    % 22  0.00000000000006
\PT{ 10.500}{ 0.001}    % 21  0.00000000000097
\PT{ 10.530}{ 0.000}    % 19 -0.00000000000025
\PT{ 10.560}{ 0.000}    % 21 -0.00000000000006
\PT{ 10.590}{ 0.000}    % 18  0.00000000000013
\PT{ 10.620}{ 0.000}    % 18 -0.00000000000017
\PT{ 10.650}{ 0.000}    % 16  0.00000000000070
\PT{ 10.680}{ 0.000}    % 10  0.00000000000061

%********************
\end{picture}}
%********************

\end{picture}

%% file: Archpict_3.tex
%Varianten:\circle* und/oder {0.04}
\def\PKT{\circle*{0.01}}
\def\PT#1#2{\put(#1,#2){\PKT}}

%%%% Hier beginnt das äußere BILD!

\unitlength1cm
\begin{picture}(13,9.5) %Breite und Höhe

% Koordinaten wie üblich: (11,1) heißt im Bild 11(cm) nach rechts und 1(cm)
% nach oben, bei Geraden ist (a,b) der Anstieg b/a und {c} die Länge,

%\put(1,1){\line(1,0){10}}
%\put(1,1){\line(0,1){8.5}}
%\put(11,1){\line(0,1){8.5}}
%\put(1,9.5){\line(1,0){10}}
\put(4.1,0.3){\hspace{1cm} {\em Diagram 3}}

%Inneres Bild
%***************************************
\put(1.5,0.9){\begin{picture}(12,8)
%***************************************
\put(0,0){\line(1,0){9.48}}
\put(0,4.68){\line(1,0){9.48}}
\put(0,3.32){\line(1,0){9.48}}
\put(4.74,0){\line(0,1){8}}

%\put(-0.5,-0.4){$X=-1.58$}
%\put(8.5,-0.4){$X=1.58$}

\put(0,-0.1){\line(0,1){0.2}}
\put(-0.2,-0.5){$-1.58$}
\put(1.74,-0.1){\line(0,1){0.2}}
\put(1.5,-0.5){$-1$}
\put(6.74,-0.1){\line(0,1){0.2}}
\put(6.75,-0.5){$1$}
\put(9.48,-0.1){\line(0,1){0.2}}
\put(9.1,-0.5){$1.58$}
\put(6.1,0.2){$X$}
\put(8.2,2.9){$\sigma=0.415$}
\put(8.2,4.8){$\sigma=0.585$}
\put(4.64,8){\line(1,0){0.2}}
\put(4.95,7.85){$1$}
\put(4.95,6.5){$\sigma$}

\PT{  0.000}{ 4.725}    % 38 -0.00000000000063
\PT{  0.030}{ 4.741}    % 33  0.00000000000047
\PT{  0.060}{ 4.757}    % 33  0.00000000000084
\PT{  0.090}{ 4.775}    % 32 -0.00000000000040
\PT{  0.120}{ 4.792}    % 33  0.00000000000093
\PT{  0.150}{ 4.810}    % 33  0.00000000000012
\PT{  0.180}{ 4.829}    % 29  0.00000000000091
\PT{  0.210}{ 4.848}    % 32 -0.00000000000054
\PT{  0.240}{ 4.868}    % 34 -0.00000000000061
\PT{  0.270}{ 4.889}    % 34  0.00000000000045
\PT{  0.300}{ 4.910}    % 32  0.00000000000077
\PT{  0.330}{ 4.932}    % 33 -0.00000000000031
\PT{  0.360}{ 4.954}    % 33  0.00000000000055
\PT{  0.390}{ 4.978}    % 31  0.00000000000028
\PT{  0.420}{ 5.002}    % 33  0.00000000000075
\PT{  0.450}{ 5.027}    % 33  0.00000000000096
\PT{  0.480}{ 5.053}    % 33 -0.00000000000047
\PT{  0.510}{ 5.081}    % 33 -0.00000000000062
\PT{  0.540}{ 5.109}    % 31 -0.00000000000056
\PT{  0.570}{ 5.140}    % 31 -0.00000000000024
\PT{  0.600}{ 5.172}    % 33 -0.00000000000059
\PT{  0.630}{ 5.206}    % 28 -0.00000000000050
\PT{  0.660}{ 5.243}    % 30  0.00000000000028
\PT{  0.690}{ 5.283}    % 33  0.00000000000085
\PT{  0.720}{ 5.327}    % 32  0.00000000000066
\PT{  0.750}{ 5.376}    % 32 -0.00000000000077
\PT{  0.780}{ 5.433}    % 32 -0.00000000000039
\PT{  0.810}{ 5.503}    % 29 -0.00000000000048
\PT{  0.840}{ 5.596}    % 30  0.00000000000058
\PT{  0.870}{ 5.790}    % 30 -0.00000000000020

\PT{  0.000}{ 7.791}    % 31 -0.00000000000082
\PT{  0.030}{ 7.755}    % 30 -0.00000000000044
\PT{  0.060}{ 7.717}    % 32 -0.00000000000020
\PT{  0.090}{ 7.677}    % 32  0.00000000000016
\PT{  0.120}{ 7.637}    % 31 -0.00000000000040
\PT{  0.150}{ 7.595}    % 31 -0.00000000000097
\PT{  0.180}{ 7.552}    % 31 -0.00000000000002
\PT{  0.210}{ 7.508}    % 30 -0.00000000000034
\PT{  0.240}{ 7.463}    % 31  0.00000000000089
\PT{  0.270}{ 7.416}    % 30 -0.00000000000015
\PT{  0.300}{ 7.369}    % 31 -0.00000000000073
\PT{  0.330}{ 7.320}    % 30 -0.00000000000087
\PT{  0.360}{ 7.270}    % 32 -0.00000000000029
\PT{  0.390}{ 7.219}    % 27 -0.00000000000064
\PT{  0.420}{ 7.167}    % 30 -0.00000000000064
\PT{  0.450}{ 7.114}    % 32 -0.00000000000026
\PT{  0.480}{ 7.059}    % 31 -0.00000000000022
\PT{  0.510}{ 7.002}    % 28  0.00000000000004
\PT{  0.540}{ 6.944}    % 30 -0.00000000000031
\PT{  0.570}{ 6.883}    % 29  0.00000000000066
\PT{  0.600}{ 6.821}    % 31  0.00000000000055
\PT{  0.630}{ 6.756}    % 31 -0.00000000000082
\PT{  0.660}{ 6.688}    % 32 -0.00000000000089
\PT{  0.690}{ 6.617}    % 32  0.00000000000031
\PT{  0.720}{ 6.541}    % 29 -0.00000000000035
\PT{  0.750}{ 6.460}    % 31 -0.00000000000026
\PT{  0.780}{ 6.370}    % 31 -0.00000000000074
\PT{  0.810}{ 6.268}    % 30 -0.00000000000060
\PT{  0.840}{ 6.142}    % 32 -0.00000000000022
\PT{  0.870}{ 5.915}    % 30  0.00000000000033
\PT{  4.485}{ 3.557}    % 29  0.00000000000026
\PT{  4.515}{ 3.679}    % 32 -0.00000000000044
\PT{  4.545}{ 3.747}    % 32  0.00000000000012
\PT{  4.575}{ 3.799}    % 30 -0.00000000000041
\PT{  4.605}{ 3.843}    % 32 -0.00000000000055
\PT{  4.635}{ 3.882}    % 32 -0.00000000000011
\PT{  4.665}{ 3.918}    % 31  0.00000000000083
\PT{  4.695}{ 3.951}    % 32  0.00000000000023
\PT{  4.725}{ 3.984}    % 30 -0.00000000000066

\PT{  4.485}{ 3.366}    % 33 -0.00000000000043
\PT{  4.515}{ 3.214}    % 30  0.00000000000096
\PT{  4.545}{ 3.115}    % 29 -0.00000000000042
\PT{  4.575}{ 3.034}    % 31 -0.00000000000090
\PT{  4.605}{ 2.962}    % 29  0.00000000000071
\PT{  4.635}{ 2.896}    % 32 -0.00000000000001
\PT{  4.665}{ 2.835}    % 31 -0.00000000000057
\PT{  4.695}{ 2.777}    % 32 -0.00000000000008
\PT{  4.725}{ 2.722}    % 29  0.00000000000013
\PT{  4.755}{ 2.669}    % 31 -0.00000000000050
\PT{  4.785}{ 2.618}    % 30 -0.00000000000056
\PT{  4.815}{ 2.569}    % 30 -0.00000000000057
\PT{  4.845}{ 2.522}    % 31  0.00000000000029
\PT{  4.875}{ 2.475}    % 31 -0.00000000000037
\PT{  4.905}{ 2.430}    % 31 -0.00000000000000
\PT{  4.935}{ 2.386}    % 24  0.00000000000003
\PT{  4.965}{ 2.343}    % 31  0.00000000000024
\PT{  4.995}{ 2.301}    % 31  0.00000000000052
\PT{  5.025}{ 2.260}    % 30 -0.00000000000053
\PT{  5.055}{ 2.220}    % 28  0.00000000000018
\PT{  5.085}{ 2.180}    % 31  0.00000000000006
\PT{  5.115}{ 2.141}    % 29 -0.00000000000058
\PT{  5.145}{ 2.103}    % 30 -0.00000000000084
\PT{  5.175}{ 2.065}    % 30  0.00000000000086
\PT{  5.205}{ 2.028}    % 31  0.00000000000020
\PT{  5.235}{ 1.991}    % 30  0.00000000000095
\PT{  5.265}{ 1.955}    % 28 -0.00000000000017
\PT{  5.295}{ 1.920}    % 30 -0.00000000000067
\PT{  5.325}{ 1.885}    % 31  0.00000000000001
\PT{  5.355}{ 1.851}    % 30  0.00000000000018
\PT{  5.385}{ 1.817}    % 31 -0.00000000000027
\PT{  5.415}{ 1.784}    % 29  0.00000000000050
\PT{  5.445}{ 1.751}    % 25 -0.00000000000004
\PT{  5.475}{ 1.718}    % 28  0.00000000000014
\PT{  5.505}{ 1.686}    % 30  0.00000000000016
\PT{  5.535}{ 1.655}    % 30 -0.00000000000005
\PT{  5.565}{ 1.624}    % 29  0.00000000000096
\PT{  5.595}{ 1.593}    % 26  0.00000000000085
\PT{  5.625}{ 1.563}    % 29  0.00000000000070
\PT{  5.655}{ 1.533}    % 29  0.00000000000036
\PT{  5.685}{ 1.503}    % 29 -0.00000000000020
\PT{  5.715}{ 1.474}    % 30 -0.00000000000033
\PT{  5.745}{ 1.445}    % 30  0.00000000000011
\PT{  5.775}{ 1.417}    % 29  0.00000000000058
\PT{  5.805}{ 1.389}    % 25 -0.00000000000067
\PT{  5.835}{ 1.361}    % 29 -0.00000000000055
\PT{  5.865}{ 1.334}    % 30  0.00000000000058
\PT{  5.895}{ 1.307}    % 30  0.00000000000018
\PT{  5.925}{ 1.280}    % 27 -0.00000000000018
\PT{  5.955}{ 1.254}    % 26 -0.00000000000001
\PT{  5.985}{ 1.228}    % 30  0.00000000000018
\PT{  6.015}{ 1.202}    % 29  0.00000000000093
\PT{  6.045}{ 1.177}    % 30 -0.00000000000006
\PT{  6.075}{ 1.152}    % 29  0.00000000000017
\PT{  6.105}{ 1.128}    % 28 -0.00000000000078
\PT{  6.135}{ 1.103}    % 26 -0.00000000000006
\PT{  6.165}{ 1.079}    % 27 -0.00000000000024
\PT{  6.195}{ 1.056}    % 30  0.00000000000061
\PT{  6.225}{ 1.032}    % 21 -0.00000000000076
\PT{  6.255}{ 1.009}    % 30 -0.00000000000015
\PT{  6.285}{ 0.987}    % 30 -0.00000000000018
\PT{  6.315}{ 0.964}    % 30 -0.00000000000032
\PT{  6.345}{ 0.942}    % 29 -0.00000000000023
\PT{  6.375}{ 0.920}    % 30 -0.00000000000044
\PT{  6.405}{ 0.899}    % 30 -0.00000000000023
\PT{  6.435}{ 0.878}    % 30  0.00000000000018
\PT{  6.465}{ 0.857}    % 29  0.00000000000023
\PT{  6.495}{ 0.836}    % 29  0.00000000000032
\PT{  6.525}{ 0.816}    % 30 -0.00000000000039
\PT{  6.555}{ 0.796}    % 30 -0.00000000000002
\PT{  6.585}{ 0.776}    % 28 -0.00000000000085
\PT{  6.615}{ 0.757}    % 30 -0.00000000000045
\PT{  6.645}{ 0.738}    % 30 -0.00000000000039
\PT{  6.675}{ 0.719}    % 29 -0.00000000000009
\PT{  6.705}{ 0.700}    % 29  0.00000000000008
\PT{  6.735}{ 0.682}    % 30  0.00000000000030
\PT{  6.765}{ 0.664}    % 28 -0.00000000000042
\PT{  6.795}{ 0.646}    % 28 -0.00000000000075
\PT{  6.825}{ 0.629}    % 28  0.00000000000072
\PT{  6.855}{ 0.612}    % 27 -0.00000000000015
\PT{  6.885}{ 0.595}    % 28  0.00000000000057
\PT{  6.915}{ 0.578}    % 28 -0.00000000000011
\PT{  6.945}{ 0.562}    % 27  0.00000000000004
\PT{  6.975}{ 0.546}    % 29 -0.00000000000029
\PT{  7.005}{ 0.530}    % 29 -0.00000000000059
\PT{  7.035}{ 0.515}    % 29  0.00000000000048
\PT{  7.065}{ 0.500}    % 28 -0.00000000000076
\PT{  7.095}{ 0.485}    % 29 -0.00000000000098
\PT{  7.125}{ 0.470}    % 27 -0.00000000000055
\PT{  7.155}{ 0.455}    % 30 -0.00000000000012
\PT{  7.185}{ 0.441}    % 29 -0.00000000000042
\PT{  7.215}{ 0.427}    % 29 -0.00000000000040
\PT{  7.245}{ 0.414}    % 29 -0.00000000000014
\PT{  7.275}{ 0.400}    % 26  0.00000000000071
\PT{  7.305}{ 0.387}    % 29  0.00000000000037
\PT{  7.335}{ 0.374}    % 29 -0.00000000000056
\PT{  7.365}{ 0.362}    % 26 -0.00000000000022
\PT{  7.395}{ 0.349}    % 28  0.00000000000016
\PT{  7.425}{ 0.337}    % 27 -0.00000000000017
\PT{  7.455}{ 0.325}    % 28 -0.00000000000073
\PT{  7.485}{ 0.314}    % 29 -0.00000000000013
\PT{  7.515}{ 0.302}    % 28  0.00000000000068
\PT{  7.545}{ 0.291}    % 24 -0.00000000000094
\PT{  7.575}{ 0.280}    % 25 -0.00000000000062
\PT{  7.605}{ 0.270}    % 21  0.00000000000064
\PT{  7.635}{ 0.259}    % 28 -0.00000000000063
\PT{  7.665}{ 0.249}    % 28  0.00000000000027
\PT{  7.695}{ 0.239}    % 25 -0.00000000000089
\PT{  7.725}{ 0.229}    % 28 -0.00000000000009
\PT{  7.755}{ 0.220}    % 28 -0.00000000000097
\PT{  7.785}{ 0.211}    % 29 -0.00000000000021
\PT{  7.815}{ 0.202}    % 29 -0.00000000000014
\PT{  7.845}{ 0.193}    % 28 -0.00000000000084
\PT{  7.875}{ 0.184}    % 27 -0.00000000000009
\PT{  7.905}{ 0.176}    % 27  0.00000000000023
\PT{  7.935}{ 0.168}    % 28  0.00000000000075
\PT{  7.965}{ 0.160}    % 28 -0.00000000000054
\PT{  7.995}{ 0.152}    % 26  0.00000000000045
\PT{  8.025}{ 0.145}    % 28  0.00000000000003
\PT{  8.055}{ 0.138}    % 28 -0.00000000000007
\PT{  8.085}{ 0.131}    % 27  0.00000000000031
\PT{  8.115}{ 0.124}    % 28  0.00000000000050
\PT{  8.145}{ 0.117}    % 26  0.00000000000042
\PT{  8.175}{ 0.111}    % 24 -0.00000000000078
\PT{  8.205}{ 0.105}    % 27 -0.00000000000046
\PT{  8.235}{ 0.099}    % 27 -0.00000000000039
\PT{  8.265}{ 0.093}    % 25 -0.00000000000087
\PT{  8.295}{ 0.088}    % 27  0.00000000000002
\PT{  8.325}{ 0.082}    % 24  0.00000000000013
\PT{  8.355}{ 0.077}    % 27  0.00000000000074
\PT{  8.385}{ 0.072}    % 27 -0.00000000000091
\PT{  8.415}{ 0.067}    % 27  0.00000000000009
\PT{  8.445}{ 0.063}    % 26 -0.00000000000081
\PT{  8.475}{ 0.058}    % 27  0.00000000000044
\PT{  8.505}{ 0.054}    % 26 -0.00000000000041
\PT{  8.535}{ 0.050}    % 25  0.00000000000045
\PT{  8.565}{ 0.046}    % 25  0.00000000000035
\PT{  8.595}{ 0.043}    % 26  0.00000000000029
\PT{  8.625}{ 0.039}    % 27  0.00000000000011
\PT{  8.655}{ 0.036}    % 26 -0.00000000000098
\PT{  8.685}{ 0.033}    % 26 -0.00000000000054
\PT{  8.715}{ 0.030}    % 26 -0.00000000000075
\PT{  8.745}{ 0.027}    % 25 -0.00000000000080
\PT{  8.775}{ 0.024}    % 25  0.00000000000083
\PT{  8.805}{ 0.022}    % 26 -0.00000000000026
\PT{  8.835}{ 0.020}    % 26  0.00000000000020
\PT{  8.865}{ 0.017}    % 25  0.00000000000036
\PT{  8.895}{ 0.015}    % 23  0.00000000000027
\PT{  8.925}{ 0.013}    % 25 -0.00000000000085
\PT{  8.955}{ 0.012}    % 25  0.00000000000008
\PT{  8.985}{ 0.010}    % 24  0.00000000000050
\PT{  9.015}{ 0.009}    % 25 -0.00000000000021
\PT{  9.045}{ 0.007}    % 22  0.00000000000009
\PT{  9.075}{ 0.006}    % 24  0.00000000000044
\PT{  9.105}{ 0.005}    % 23 -0.00000000000091
\PT{  9.135}{ 0.004}    % 23 -0.00000000000017
\PT{  9.165}{ 0.003}    % 23  0.00000000000081
\PT{  9.195}{ 0.003}    % 22  0.00000000000079
\PT{  9.225}{ 0.002}    % 19 -0.00000000000019
\PT{  9.255}{ 0.001}    % 18  0.00000000000013
\PT{  9.285}{ 0.001}    % 18 -0.00000000000049
\PT{  9.315}{ 0.001}    % 18 -0.00000000000009
\PT{  9.345}{ 0.000}    % 20 -0.00000000000013
\PT{  9.375}{ 0.000}    % 18 -0.00000000000017
\PT{  9.405}{ 0.000}    % 17 -0.00000000000097
\PT{  9.435}{ 0.000}    % 16 -0.00000000000039
\PT{  9.465}{ 0.000}    % 13  0.00000000000003

%********************
\end{picture}}
%********************

\end{picture}